\documentclass[12pt]{amsart}

\usepackage{amsfonts}
\usepackage{enumerate}
\usepackage{verbatim}
\usepackage{graphicx}
\usepackage{subfig}
\usepackage{wrapfig}
\usepackage{array}
\usepackage{booktabs}
\usepackage{bbold}
\usepackage[T1]{fontenc}
\usepackage{tikz}
\usetikzlibrary{matrix}

\newcommand{\ind}{$\quad$}
\newcommand{\odpar}{\vspace{0.5cm}}
\newcommand{\CS}{C(\mathcal{S})}

\newcommand\numberthis{\addtocounter{equation}{1}\tag{\theequation}}

\newcommand{\N}{\mathbb{N}}
\newcommand{\R}{\mathbb{R}}
\newcommand{\Sp}{\mathbb{S}}
\newcommand{\Sd}{\mathbb{S}^d}
\newcommand{\Z}{\mathbb{Z}}
\newcommand{\Cal}{\mathcal}
\newcommand{\calD}{{\mathcal {D}}}

\newcommand{\supp}{\operatorname{supp}}
\newcommand{\spa}{\operatorname{span}}

\newcommand{\ch}{\mathbf 1}

\newcommand{\ve}{\varepsilon}

\newcommand{\lan}{\langle}
\newcommand{\ran}{\rangle}
\newcommand{\ints}{\int_{\Sd}}
\newcommand{\no}{\nonumber\\}
\newcommand{\calJ}{{\mathcal {J}}}
\newcommand{\jn}{j_{n}}

\def\Koniec{\hbox to\hsize{\hfil$\diamond$}}

\def\1{\mathbf 1}
\def\K{\Cal{K}}
\def\E{\mathcal{E}}

\DeclareRobustCommand{\Biess}[0]{B^{s}_{2,\infty}(\Sp^d)} 

\DeclareRobustCommand{\nk}[1]{\left\|{#1}\right\|_{\Cal{K}}} 
\DeclareRobustCommand{\dS}[1]{d\sigma_d{(#1)}} 
\DeclareRobustCommand{\kom}[1]{\tag*{ $\langle$ #1 $\rangle$}}
\newcommand\num{\addtocounter{equation}{1}\tag{\theequation}} 
\DeclareRobustCommand{\ll}[1]{\left\|{#1}\right\|_2} 
\DeclareRobustCommand{\CE}[0]{C_{_E}} 
\DeclareRobustCommand{\D}[0]{D_{_N}} 
\DeclareRobustCommand{\M}[0]{M_{_N}} 
\DeclareRobustCommand{\nb}[1]{\left\|{#1}\right\|_{s,2}} 

\newtheorem{theorem}{Theorem}[section]
\newtheorem{corollary}[theorem]{Corollary}
\newtheorem{lemma}[theorem]{Lemma}

\theoremstyle{definition}
\newtheorem{definition}{Definition}[section]
\theoremstyle{remark}

\numberwithin{equation}{section}

\begin{document}

\title[]{Multiresolution analysis and adaptive estimation on
a sphere using stereographic wavelets}

\author{Bogdan \'Cmiel}
\address{Faculty of Applied Mathematics, AGH University of Science and Technology, Al. Mickiewicza 30, 30-059 Cracow, Poland.}
\email{cmielbog@gmail.com}

\author{Karol Dziedziul}
\address{
Faculty of Applied Mathematics,
Gda\'nsk University of Technology,
ul. G. Narutowicza 11/12,
80-952 Gda\'nsk, Poland}

\email{karol.dziedziul@pg.edu.pl}

\author{Natalia Jarz\k{e}bkowska}
\address{
Faculty of Applied Mathematics,
Gda\'nsk University of Technology,
ul. G. Narutowicza 11/12,
80-952 Gda\'nsk, Poland}




\date{\today}

\begin{abstract} 
We construct an adaptive estimator of a density function on $d$ dimensional unit sphere $\Sp^d$
($d \geq 2 $), using a new type of spherical frames. The frames, or as we call them, stereografic wavelets are obtained by transforming a wavelet system, namely Daubechies, using some stereographic operators. We prove that our estimator achieves an optimal rate of convergence on some Besov type class of functions by adapting to unknown smoothness. Our new construction of stereografic wavelet system gives us a multiresolution approximation of $L^2(\Sp^d)$ which can be used in many approximation and estimation problems.
In this paper we also demonstrate how to implement the density estimator in $\Sp^2$ and we present a finite
sample behavior of that estimator in a numerical experiment.

\end{abstract}

\maketitle

\section{Introduction}

In this paper, we consider an adaptive estimator of a density function on the $d$-dimensional unit sphere $\Sp^d$, $d \geq 2 $ using a new type of Parseval frame. To construct the estimator we create a new stereografic wavelet system which gives us a multiresolution approximation of $L^2(\Sp^d)$. Since our construction uses a standard wavelet system (namely Daubechies) and some stereografic operators one can only make some modifications of existing algorithms in $\mathbb{R}^d$, which is relatively easy, to enjoy the benefits of mutiresolution analysis on a sphere and solve many approximation and estimation problems. 

Let us start from the definition
\begin{definition}\label{def:estymator}
Let $\{K_j: j\geq j_0\}$ be a family of measurable functions (called kernels) $K_j:\Sd\times \Sd \to \R$. 
Let  $X_1,\ldots X_n$ be i.i.d. with density function $f$ on $\Sd$ with respect to Lebesgue'a measure. For $j\geq j_0$ we define  an estimator of $f$ 
\[
f_n(j)(x)=\frac{1}{n}\sum_{i=1}^n K_j(x,X_i).
\]
 \end{definition}
\noindent
We denote the balls of density functions in Besov spaces on sphere
\[
\Sigma(s,\tilde{B})=\{f\in \Biess: \int_{\Sd} f(x)\dS{x}=1,\; f\geq 0,\; \|f\|_{s,2}\leq \tilde{B} \}.
\]
Since we want to obtain an adaptive estimator, we want to construct kernels $K_j$ on a sphere for which we have an optimal rate of estimation. Namely,
\begin{theorem}\label{twierdzenie0}
Let $d/2< r<R$ and let   $X_1,\ldots X_n$ be i.i.d.  with density function $f\in \Biess$, where $s$ is unknown and ${r\leq s\leq R}$.
 Then there is a family of kernels $\{K_j: j\geq j_0\}$ such that for any $U>0$ there are constants $c=c(r,R,U)$ and $C=C(U)$ such that for all $s, n$ and $\tilde{B}>1$ we have 
\begin{align}\nonumber
\sup_{f\in \Sigma(s,\tilde{B}),\|f\|_\infty\leq U} \mathbb{E}\|f_n(j_n)-f\|_2^2\leq c \tilde{B}^{2d/(2s+d)}n^{-2s/(2s+d)}, 
\end{align}
where 
\begin{align}\nonumber
\jn=\min\left\{j\in [j_{min},j_{max}]: \, \forall_{l, \ j < l \leq j_{max}} \;\ll{ f_n(j)-f_n(l) }^2 \leq C \frac{2^{ld}}{n}\right\} 
\end{align}
and
\begin{align}\nonumber
j_{min}&=\left\lfloor \frac{\log_2n}{2R+d} \right\rfloor,\qquad
j_{max}=\left\lceil \frac{\log_2n}{2r+d} \right\rceil.
\end{align}
\end{theorem}
In the above theorem the smoothness parameter $s$ is unknown but for choosing the resolution level we use a lower bound $r$ and an upper bound $R$. Let us discuss some consequences of choosing different values for $r$ and $R$. It seems that it is good idea to take $r$ as small as possible and $R$ as big as possible to consider a very wide range for the unknown smoothness. The first part of that is true since there are no serious consequences of taking small $r$. Unfortunately if we take a big value for $R$, then we need to use in our construction some very smooth wavelets (smoother than $R$). The smoother the wavelets are, the bigger support they have and if one scales them to a smaller area, then they change values very rapidly. In the asymptotic point of view this is not a problem but for fixed $n$ the estimator loses its efficiency if the value $R$ is too big. The same problem we can observe in case of a wavelet estimation on $\mathbb{R}$. 

It is well-known (see Hall, Kerkyacharian and Picard (1998) \cite{Hall1} Theorem 4.1) that on the real line if one
considers wavelets estimators with a block thresholding procedure, one attains minimax
rate of convergence without extraneous logarithmic factors for $B^s_{2,\infty}$ Besov spaces and
$L^2$-loss, i.e., $n^{-2s/(1+2s)}$. Similar result was given in \cite{BN}. We follow the arguments presented there.

The problem of estimating nonparametrically a density on the $d$-dimensional unit
sphere $\Sd$ over Besov classes is not new (see Baldi, Kerkyacharian, Marinucci and Picard
(2009) \cite{Baldi} for a direct setting and for an indirect setting see Kerkyacharian, Pham Ngoc
and Picard (2011), \cite{kerkyacharian1}). 
In particular, in Baldi, Kerkyacharian, Marinucci and Picard (2009), the
authors had already dealt with the considered problem in a more general
framework, namely, by considering $B^s_{q,r}$ Besov spaces. 
They constructed an adaptive estimator based on a set of spherical wavelets, named
needlets, with a hard thresholding procedure. They obtained minimax rates of convergence for
$B^s_{q,r}$ Besov spaces, $L^p$-loss and sup-norm loss up to a logarithmic factor.
This deep approach was continued in \cite{Dur} but in regression case. Moreover the rates are
without the logarithmic factor.
We obtain the minimax rate of convergence
for the $L^2$-loss, i.e., $n^{-2s/(2s+d)}$ without the logarithmic factor that one usually gets with adaptive methods of estimating density function on the sphere. We want to emphasize that section 2 can be rewritten  for a compact smooth manifold.  So if one can construct a family of kernels such that they satisfy conditions of Theorem \ref{tw:TalagrandMOJ},  one obtains a method of minimax rate of convergence
for the $L^2$-loss, i.e., $n^{-2s/(2s+d)}$ without the logarithmic factor on the manifold $M$. The first step is done in \cite{BDK} i.e., a smooth orthogonal decomposition of identity in $L^2(M)$ is constructed.

In section 2 we formulate  general  conditions on $K_j$ 
which  guarantee Talagrand's inequality in Bousquet's version. This inequality is a key argument to prove Theorem \ref{twierdzenie0}.

In section 3 we construct kernels $K_j$ using a new type of frames on sphere which gives an optimal rate of the estimation. 
The new frame, called stereographic wavelets, inherits all properties of the classical  multivariate Daubechies wavelets and we know that such basis is an excellent tool in the process of approximation and estimation mentioned in \cite[Theorem 2]{BN} or \cite[Theorem 4.1]{Hall1}.

In great amount of literature, tight frames (needlets) are used as a tool in approximation as well as in estimation of densities (see \cite{Dai}, \cite{geller}, \cite{kerkyacharian}, \cite{NP2}, \cite{NPW}).
Unfortunately this approach does not give an optimal rate of estimation (see \cite{Kueh}).
The new frames on sphere, introduced by Bownik M., Dziedziul K. in \cite{BD} give a construction of  $K_j$ such that we achieve the optimal rate of convergence on Besov spaces via adaptive estimation.  
 The method of constructing Parseval frame on sphere consists of two steps. In first step we obtain a localized wavelet system on sphere by transforming Daubechies wavelet system on $[1-\ve,1+\ve]^d$ using two stereographic operators. Next we create a Parseval frame by applying P. Auscher, G. Weiss, M. V. Wickerhauser (AWW) operator (\cite{AWW}) on sphere (see \cite{BD}). 
Consequently we create a multivariate approximation on $L^2(\Sp^d)$, i.e.,
\[
V_{j_0}(\Sd)\subset  V_{j_0+1}(\Sd)\subset \cdots \subset L^2(\Sp^d),
\]
for all $j\geq j_0$
\[
\dim V_j(\Sd)<\infty,
\]
and $\bigcup_{j\geq j_0} V_j(\Sd)$ is dense in $L^2(\Sd)$.
The functions $K_j(x,y)$, $x,y\in \Sp^d$, $j\geq j_0$ in the main theorem are kernels of the orthogonal projection 
\[
K_j: L^2(\Sp^d) \to V_j(\Sd),
\]
\[
K_j (f)(x)=\int_{\Sp^d} f(y) K_j(x,y) d\sigma_d(y),
\]
where $\sigma_d$ is Lebesgue measure on $\Sd$.

In section 4 we present a technical version of Theorem \ref{twierdzenie0}. In section 5 we show a numerical example of such estimation for $\Sp^2$ (classical sphere). All the proofs are given in the appendix A.

Since the coefficients of the frame give us characterization of  Besov spaces  $\Biess$, it is possible to 
use our approach from earlier papers to estimate the smoothness of density function or to construct a smoothness test (see  \cite{CD}, \cite{CDW} and \cite{kucharska}), but this is not the aim of this paper.

\section{ Talagrand's inequality}
\label{cha:Talagrand}

In this section we present Talagrand's inequality (see \cite{Boucheron}, \cite{gn} and Theorem 3.3.9 (Upper tail of Talagrand's inequality, Bousquet's version \cite{Gine}) and its consequences. Let us cite from  \cite{gn}: in the special case "Talagrand's inequality becomes exactly the Bernstein and Prohorov inequalities. Clearly then,  Talagrand's inequality is essentially a best possible exponential bound for the empirical process."
We start  our consideration with general type of kernels $K_j(x,y)$, $x,y\in \Sp^d$. 
In the next section we focus our attention on  kernels which arise from the Parseval frame.

\begin{theorem}\label{tw:talagrand}
Let $X,X_1, \ldots, X_n$ be i.i.d. random variables with law  $\mu$ on a measurable space $(M,\Cal{M}).$ Let $\Cal{K}$ be a countable class of real measurable functions on  $M$, uniformly bounded by a constant $U$ and $\mu-$ centered, i.e.,
\begin{eqnarray}\label{zal:1_scen_tw:talagrand}
\mu(k):=\int kd\mu=0, \quad k\in\Cal{K}.
\end{eqnarray}
For $H:\Cal{K} \rightarrow \R$ define
\begin{eqnarray}\label{def:normaF_tw:talagrand}
\nk{H} &=& \sup_{k\in\Cal{K}} |H(k)|.
\end{eqnarray}
Let $\omega$ be a positive number such that 
\begin{eqnarray}\label{zal:2_omega_tw:talagrand}
\omega^2&\geq& \sup_{k\in\Cal{K}} \mathbb{E}k^2(X)
\end{eqnarray}
and
\begin{eqnarray}\label{zal:3_V_tw:talagrand}
V&:=&n\omega^2 + 2U\mathbb{E}\nk{\sum_{i=1}^n k(X_i)}.
\end{eqnarray}
Then for every $x\geq 0$ and $n\in\N$
\begin{eqnarray}\label{eq:talagrand}
P\left\{  \nk{\sum_{i=1}^n k(X_i)} \geq \mathbb{E}\nk{\sum_{i=1}^n k(X_i)} + \sqrt{2Vx} +Ux/3    \right\}&\leq&2e^{-x}.
\end{eqnarray}

\end{theorem}


Following \cite{gn} we adapt this theorem to our situation.
Let $X_1, \ldots X_n$ be i.i.d. random variables with density $f$ on $\Sd$ with respect to 
 Lebesgue measure $\sigma_d$ on  $\Sd$ (the surface measure),
$$
\Sd = \left\{ x \in \R^{d+1}: ||x||=\sqrt{x_1^2+\ldots + x_{d+1}^2}=1  \right\}\subset \R^{d+1}.
$$ 
Let us define
$$
d\mu=fd\sigma_d.
$$
We assume that
$$
\|f\|_\infty<\infty.
$$
Since  $L^2(\Sd)$ is a separable set, then there is a countable subset $B_0$ of the unit ball  $B\subset L^2(\Sd)$ such that for all $\beta\in L^2(\Sd)$ 
\begin{align}\label{eq:Hahn_Banach_tw:talagrand}
\ll{\beta} &= \left(\int_{\Sd} |\beta(t)|^2 \dS{t}\right)^{1/2} = \sup_{g\in B_0} \left|  \ints \beta(t)g(t) \dS{t}\right|=\sup_{g\in B_0} \left| \lan\beta,g \ran\right|
\end{align}
We assume that a family of symmetric kernels  $K_j(\cdot,\cdot)$, $j\geq j_0$ i.e.
$K_j(x,y)=K_j(y,x)$ for all $x,y,j$,  satisfies the following three conditions:
\begin{equation}\label{war1}
\forall_{j\geq j_0} \quad \sup_{x,y\in \Sd}|K_j(x,y)|< \infty,
\end{equation}
\begin{equation}\label{war2}
\exists_{D>0} \forall_{y\in \Sd} \forall_{j\geq j_0} \quad \ints K_j^2(x,y)\dS{x}\leq D 2^{jd},
\end{equation}
\begin{equation}\label{war3}
\exists_{C_\K>0}\forall_{j\geq j_0}\forall_{g\in L^2(\Sd)}\quad \ll{K_j(g)}\leq C_{\K} \ll{g},
\end{equation}
where 
\begin{align}\nonumber
K_j(g)(t)=\ints g(x)K_j(x,t)\dS{x},
\end{align}
for $g\in L^1(\Sd)$. Note that 
\begin{align}\label{eq:K_j}
K_j(f)(t)=\mathbb{E}[K_j(t,X)].
\end{align}
For simplicity we will assume further that $C_{\K}=1.$
A classical estimator of density $f$ is given by
\[
f_n(j)(x)=\frac{1}{n}\sum_{i=1}^n K_j(x,X_i).
\]
From  (\ref{war2}) we obtain the following lemma
\begin{lemma}\label{lem:SzacowanieWariancji}
Let $X_1, \ldots X_n$ be i.i.d. with common density $f$ on $\Sd$ with respect to Lebesgue measure. 
Let  symmetric kernels $K_j(\cdot,\cdot)$ satisfy (\ref{war1}), (\ref{war2}), (\ref{war3}). Then, there exists $D>0$ such that
\begin{eqnarray}\label{lSW1}
\forall_{j\geq j_0} \quad \mathbb{E}\ll{f_n(j)-\mathbb{E}f_n(j)}^2\leq D\frac{2^{jd}}{n},
\end{eqnarray}
and
\begin{eqnarray}\label{lSW2}
\forall_{j\geq j_0} \quad \mathbb{E}\ll{\sum_{i=1}^n \Big(K_j(\cdot,X_i)- \mathbb{E}K_j(\cdot,X_i)\Big)}\leq  \sqrt{Dn2^{jd}}.
\end{eqnarray}
{\emph{[proof in the appendix \ref{A1}]}}
\end{lemma}
\noindent
For  $j_0\in \Z$ and  $j\geq j_0$ we define the following family of kernels
\begin{align}\label{eqdef:rodzinaK}
\Cal{K}=\Cal{K}_j &:= \left\{k_g = \ints g(t) K_j(t,\cdot)\dS{t} - \ints g(t)K_j(f)(t)\dS{t} : g\in B_0  \right\}.
\end{align}

\begin{lemma}\label{assumption}
Let $X_1, \ldots X_n$ be i.i.d. with common density $f$ on $\Sd$ with respect to Lebesgue measure.
Let  symmetric kernels $K_j(\cdot,\cdot)$ satisfy (\ref{war1}), (\ref{war2}), (\ref{war3}).
Then $\Cal{K}$ satisfies the assumptions of Theorem \ref{tw:talagrand} i.e.,
for all $k_g\in \Cal{K}$
\begin{eqnarray}\label{eqdef:U_K}
\|k_g\|_\infty&\leq& \sqrt{D}2^{jd/2} + \|f\|_\infty^{1/2} =: U_{\Cal{K}_j}<\infty.
\end{eqnarray}
\begin{eqnarray}\label{eq:scentrowane talag}
\mu(k_g) = \ints k_g(x)f(x)\dS{x}=0.
\end{eqnarray}
\begin{eqnarray}\label{szacowanie}
 \mathbb{E}\big[k_g(X)\big]^2  
\leq \|f \|_\infty =: \omega_{\Cal{K}}^2.
\end{eqnarray}
{\emph{[proof in the appendix \ref{AA}]}}
\end{lemma}
To transform the thesis of Theorem \ref{tw:talagrand} in    a case of density function estimation note that if we define $\mu_n = \frac{1}{n}\sum_{i=1}^n \delta_{X_i}$, then for every $k_g\in \Cal{K}$
\begin{eqnarray}
\mu_n(k_g) &=&  \frac{1}{n}\sum_{i=1}^n k_g(X_i)\nonumber\\
&=& \ints g(t)\left( \frac{1}{n}\sum_{i=1}^n K_j(t,X_i) -  K_j(f)(t) \right)\dS{t}\label{eq:sum_kg_Xi}\\
&=& \ints g(t) \Big(f_n(j)(t) - \mathbb{E}f_n(j)(t)\Big)\dS{t}\label{eq:mu_n_kg}.
\end{eqnarray}
Hence taking $H=\mu_n-\mu$ in Theorem \ref{tw:talagrand} we have
\begin{align*}
\nk{\mu_n-\mu} &=\; \sup_{k_g\in \Cal{K}} \Big|  (\mu_n-\mu)k_g \Big|\kom{ \eqref{def:normaF_tw:talagrand}}\\
&= \sup_{g\in B_0} \left|   \ints g(t) \Big(f_n(j)(t) - \mathbb{E}f_n(j)(t)\Big)\dS{t}\right|\kom{ \eqref{eq:scentrowane talag} and \eqref{eq:mu_n_kg}}\\
&=\ll{f_n(j)-\mathbb{E}f_n(j)}=\frac{1}{n}\ll{\sum_{i=1}^n \Big(K_j(\cdot,X_i)- \mathbb{E}K_j(\cdot,X_i)\Big)}. \kom{ \eqref{eq:Hahn_Banach_tw:talagrand}}\\
\end{align*}

\noindent
Note that by \eqref{eq:sum_kg_Xi}
\begin{align}
\nk{\sum_{i=1}^n k_g(X_i)} 
&= n\ll{f_n(j)-\mathbb{E}f_n(j)}.\label{eq:normaK_sumy_kg_talagrand}
\end{align}

Now taking into account the above result we can formulate the following one.
\begin{corollary}\label{wn:wniosekTalagrand}
Let $X_1, \ldots X_n$ be i.i.d. with common density $f$ on $\Sd$ with respect to Lebesgue measure. 
Let  symmetric kernels $K_j(\cdot,\cdot)$ satisfy (\ref{war1}), (\ref{war2}), (\ref{war3}).
For the family $\Cal{K}=\Cal{K}_j$  \eqref{eqdef:rodzinaK} we have the following inequality (Talagrand's inequality \eqref{eq:talagrand} from Theorem \ref{tw:talagrand}) 
\begin{align}\label{eq:talagrandS}
&\forall_{x\geq 0}\forall_{n\in\N}\forall_{j\geq j_0}\no 
&\qquad P\left\{ n\ll{f_n(j)-\mathbb{E}f_n(j)} \geq n\mathbb{E}\ll{f_n(j)-\mathbb{E}f_n(j)}+ \sqrt{2Vx} +U_{\Cal{K}_j}x/3    \right\}\leq 2e^{-x},\no 
\end{align}
where
\begin{align}\label{eqdef:V_K}
V
&=n\omega_{\Cal{K}}^2 + 2nU_{\Cal{K}_j}\mathbb{E}\ll{f_n(j)-\mathbb{E}f_n(j)}
\end{align}
and
\begin{align}
\omega_{\Cal{K}}^2&=\|f\|_\infty\label{eqdef:omegaK},\\
U_{\Cal{K}_j}&= \sqrt{D}2^{jd/2} + \|f\|_\infty^{1/2}.\label{eqdef:UK}
\end{align}
\end{corollary}

\odpar
We want to transform \eqref{eq:talagrandS} into a formula which will be convenient in our later 
calculation. 
Note that
\begin{align*}\num\label{eq:2Vx}
\sqrt{2Vx}
&=\sqrt{2nx\omega_{\Cal{K}}^2 + 4xU_{\Cal{K}_j} n\mathbb{E}\ll{f_n(j)-\mathbb{E}f_n(j)} }\\
&\leq\sqrt{2nx\omega_{\Cal{K}}^2} + 2\sqrt{xU_{\Cal{K}_j} n\mathbb{E}\ll{f_n(j)-\mathbb{E}f_n(j)} }\kom{ $\sqrt{a+b}\leq \sqrt{a}+\sqrt{b} $}\\
&\leq\sqrt{2nx\omega_{\Cal{K}}^2} + xU_{\Cal{K}_j} + n\mathbb{E}\ll{f_n(j)-\mathbb{E}f_n(j)}.  \kom{ $\sqrt{ab}\leq \frac{a+b}{2} $} \\
\end{align*}
By \eqref{eq:2Vx}, \eqref{eqdef:omegaK}, \eqref{eqdef:UK}
\begin{align*}
RHS 
&:=n\mathbb{E}\ll{f_n(j)-\mathbb{E}f_n(j)}+ \sqrt{2Vx} +U_{\Cal{K}_j}x/3  \\
&\leq 2n\mathbb{E}\ll{f_n(j)-\mathbb{E}f_n(j)} + \frac{4}{3} (\sqrt{D2^{jd}} + \|f\|_\infty^{1/2})x+\sqrt{2xn\|f\|_\infty}
\end{align*}
Note that by Lemma \ref{lem:SzacowanieWariancji}
\begin{align*}
n\mathbb{E}\ll{f_n(j)-\mathbb{E}f_n(j)}&= \mathbb{E}\ll{\sum_{i=1}^n (K_j(\cdot,X_i) - \mathbb{E}K_j(\cdot,X_i))} 
\leq\sqrt{Dn2^{jd}}.
\end{align*}
Consequently, for $x=2^{jd}$ we obtain
\begin{align}\label{eq:AB_talagrand_cd}
RHS&\leq n  \frac{2^{jd/2}}{\sqrt{n}}\left(2\sqrt{D} + \frac{4}{3}\sqrt{D} \frac{2^{jd}}{\sqrt{n}} + \frac{4}{3}\|f\|_\infty^{1/2}\frac{2^{jd/2}}{\sqrt{n}}+\sqrt{2\|f\|_\infty}\right).
\end{align}
It follows from Lemma  \ref{lem:SzacowanieWariancji}, that we need to assume that the relation between level $j$ and the size of sample $n$ is such that
\[
\frac{2^{jd}}{n} \to 0
\] as $j,n \to \infty$.
This condition is
justified to guarantee balance between stochastic and deterministic error.
 So through this paper we assume that  there is $\CE\geq 0$ such that
we work in the range of parameter $j$ and $n$ such that
\begin{equation}\label{CE}
\frac{2^{jd}}{n}\leq \CE.
\end{equation}
Consequently
\begin{align}
RHS&\leq n \frac{2^{jd/2}}{\sqrt{n}}M\sqrt{1\vee\|f\|_\infty}, \label{eq:RHSostatni} 
\end{align}
where $M=2(\sqrt{D}\left(2 + \frac{4}{3} \CE\right)\ \vee \ (\frac{4}{3}\CE+\sqrt{2})).$\\
\\
Finally from \eqref{eq:RHSostatni} and Corollary \ref{wn:wniosekTalagrand} we get the main estimation. 
If we denote by $\mathcal{L}(\Sd)$ a $\sigma-$ algebra of  Lebesgue sets contained in $\Sd$, then the following theorem is true. 
\begin{theorem}\label{tw:TalagrandMOJ}
Let  $K_j(\cdot,\cdot), j\geq j_0$ be a family of real symmetric measurable function with respect to $\mathcal{L}(\Sd)\times \mathcal{L}(\Sd),$ satisfying:
\begin{equation}\label{warTW:jadro_1_ogr}
\forall_{j\geq j_0} \quad \sup_{x,y\in \Sd}|K_j(x,y)|< \infty, 
\end{equation}
\begin{equation}\label{warTW:jadro_2_calkazkwadratu}
\exists_{D>0} \forall_{y\in \Sd} \forall_{j\geq j_0} \quad \ints K_j^2(x,y)\dS{x}\leq D 2^{jd}.
\end{equation}
\begin{equation}\label{warTW:jadro_3_ograniczone}
\forall_{j\geq j_0}\forall_{g\in L^2(\Sd)}\quad \ll{K_j(g)}\leq \ll{g}.
\end{equation}
Let $X_1, \ldots X_n$ be i.i.d. with common density $f\in L^\infty(\Sd)$. Let
\[
f_n(j)(x)=\frac{1}{n}\sum_{i=1}^n K_j(x,X_i).
\]
Then for $j\geq j_0$ and $n\in N$ such that $\frac{2^{jd}}{n}\leq \CE$ we have
\begin{align*}
\quad P\left\{ \ll{f_n(j)-\mathbb{E}f_n(j)} \geq  \frac{2^{jd/2}}{\sqrt{n}}M\sqrt{1\vee\|f\|_\infty}
\right\}&\leq 2e^{-2^{jd}},
\end{align*}
where 
$$
M=2\max\left\{\sqrt{D}\left(2 + \frac{4}{3} \CE\right),\frac{4}{3}\CE+\sqrt{2}\right\}.
$$
\end{theorem}

\section{Besov spaces and  stereographic wavelets on sphere}

\begin{figure}[ht]
\centerline{\scalebox{1}
       {\includegraphics[trim = 80mm 40mm 80mm 25mm, scale=0.28]{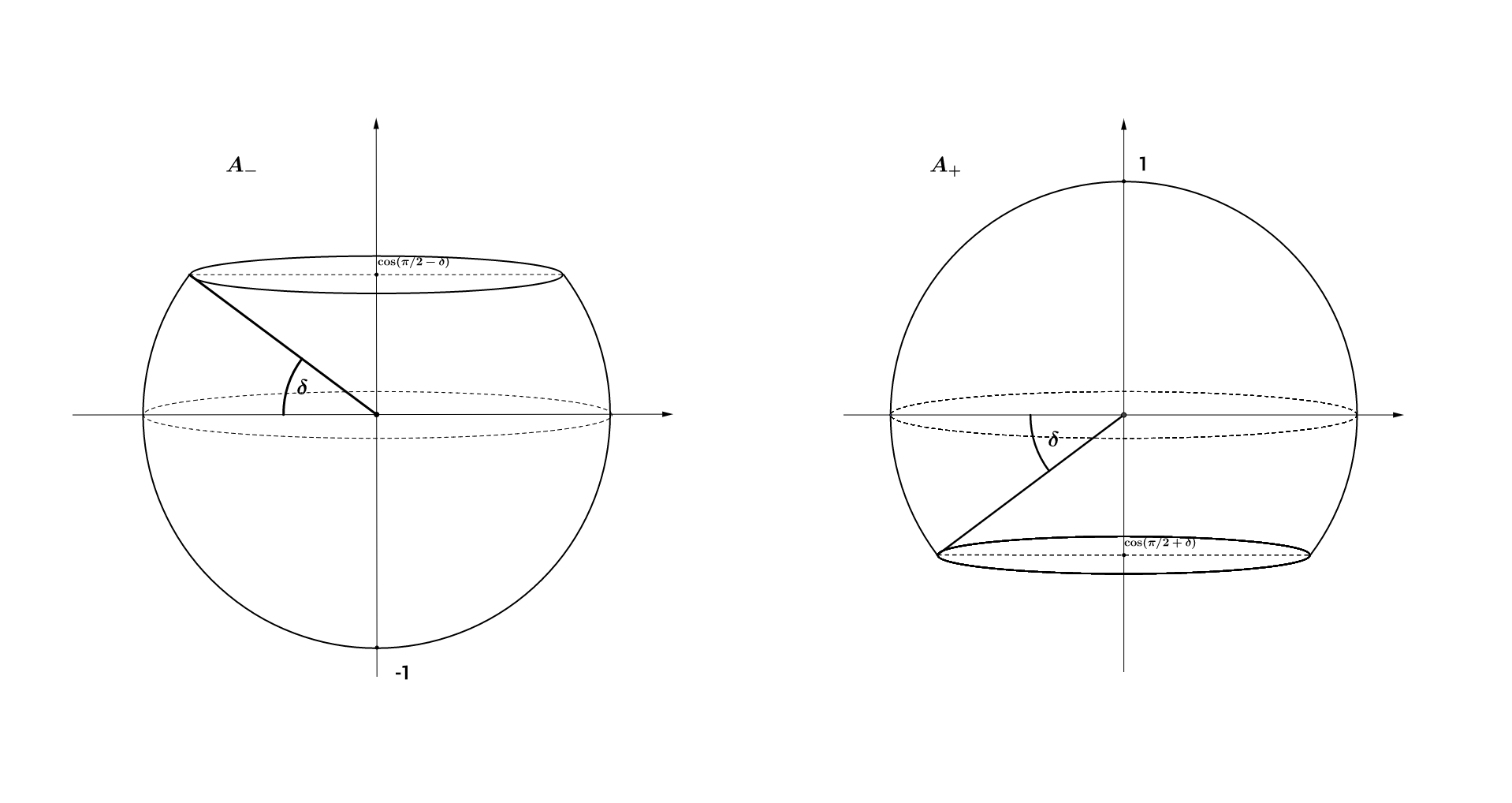}}}
       \caption{Patches $A_-$ i $A_+$ on $\Sd$}
       \label{rys:laty}
\end{figure}

\ind In this section we construct a Parseval frame on $\Sd, d\geq 2$ using Bownik-Dziedziul construction (see \cite{BD}).
For any fixed angle $0<\delta<\pi/2$ let us decompose the sphere onto two Patches
 $A_-$ and $A_+$, depending on the angle (see Figure \ref{rys:laty}), where
\begin{align*}
 A_-&=\{x \in \Sp^d: x_{d+1}\leq \cos(\pi/2-\delta) \},\\
 A_+&=\{x \in \Sp^d: x_{d+1}\geq\cos(\pi/2+\delta) \}.
\end{align*}

\begin{figure}[ht]
\centerline{\scalebox{1.07}
       {\includegraphics[scale=0.3]{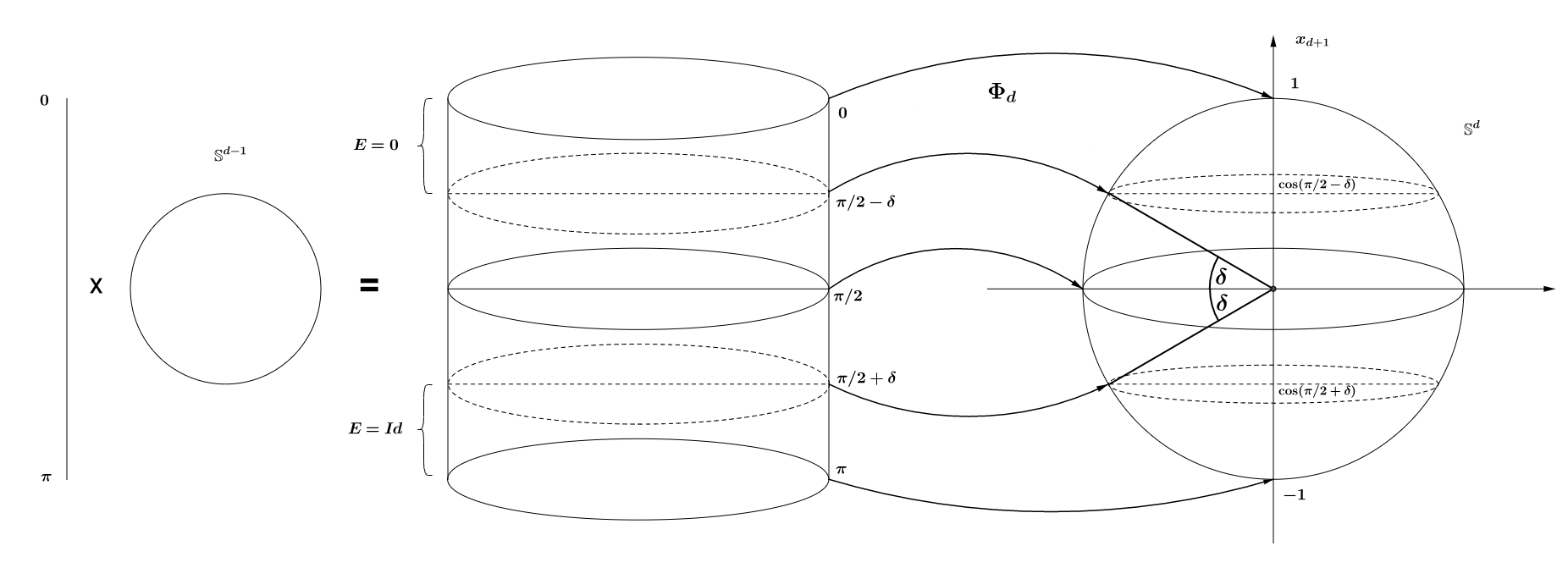}}}
       \caption{Function $\Phi_d$ and operator $E$}
       \label{rys:walec}
\end{figure}

Consider a natural parametrization of the sphere
$$
\Phi_d:[0,\pi]\times \Sp^{d-1}\rightarrow \Sp^d,\qquad
\Phi_d(\theta,\xi)=(\xi\sin\theta,\cos\theta),\quad (\theta,\xi)\in[0,\pi]\times \Sp^{d-1}.
$$
The function $\Phi_d:(0,\pi)\times \Sp^{d-1}\rightarrow \Sp^d\setminus\{\textbf{1}^d,-\textbf{1}^d\}$ is a diffeomorphism, where 
\[
\textbf{1}^d=(0,\ldots,0,1)\in\Sd
\quad  \text{is the "North Pole", see Figure \ref{rys:walec}.}
\]
Then  $f\in L^1(\Sd,d\sigma_d)$ and we have (see \cite[(1.5.4)]{xu})
\begin{align*}
\ints f(u)\dS{u} = \int_{\mathbb{S}^{d-1}}\int_0^\pi f \circ \Phi_d(\theta,\xi)(\sin\theta)^{d-1}d\theta d\sigma_{d-1}(\xi).
\end{align*}
It makes sense to introduce the notation $g(\theta,\xi)=g(\Phi_d(\theta,\xi))$. Let us take some real-valued, smooth function $s\in C^{\infty}(\R)$ such that
$$
\supp s\subset [-\delta,\infty),\qquad s^2(t)+s^2(-t)=1,\quad t\in\R.
$$
Now we can define Auscher-Weiss-Wickerhouser (AWW) operator $E=E_{\delta,s}$, pointwise for every $g:\Sp^d\rightarrow\R$
\begin{equation}\label{def:AWW}
E(g)(\theta,\xi)=\left\{ \begin{array}{l}
g(\theta,\xi), \qquad\theta>\pi/2+\delta\\
\\
s^2(\theta-\pi/2)g(\theta,\xi) + s(\theta-\pi/2)s(\pi/2-\theta)g(\pi-\theta,\xi)\\
\\
0,\qquad \theta<\pi/2-\delta,\\
\end{array}\right.
\end{equation}
where $\xi\in\Sp^{d-1}$ (see (3.5) and (3.6) in \cite{BD}). 

\begin{lemma}{\emph{\cite[Lemma~3.3]{BD}}}
$E:L^2(\Sp^d,d\sigma_d)\rightarrow L^2(\Sp^d,d\sigma_d)$ is an orthogonal projection.
\end{lemma}
Since we consider two patches  $A_-$ and $A_+$ we have very simple decomposition of identity operator $I$
(orthogonal partition of unity), see \cite{BD} Theorem 1.1
 $$
 E^+=Id-E \quad\text{ and }\quad E^-=E.
 $$
\begin{definition}
We say that an operator $P$ is localized on an open set $U$, i.e., for any $f:\mathbb S^d \to \R$ we have
\[
P f(x) = 0 \qquad\text{for }x\in \mathbb S^d \setminus U.
\]
\end{definition}
Now we are ready to reformulate \cite[Theorem~5.1]{BD}. We recall that Sobolev space \cite{He}.
\begin{definition}\label{sob}
For $k\in \N$ and $f: \Sd \to \R$ we denote by $\nabla^k f(x)$, $x\in \Sd$, the covariant derivative of $f$ of order $k$ in some local chart. We let $|\nabla^k f|$ to be its norm (which is independent of a choice of chart). 
Given $1\le p< \infty$ we define the norm 
\[
||f||_{W^r_p} = \sum_{k=0}^r \bigg( \int_{\Sd} |\nabla^k f(x)|^p d\sigma_d(x) \bigg)^{1/p}<\infty.
\]
The  Sobolev space $W^r_p(\Sd)$ is the completion of $C^r(\Sd)$  with respect to the norm $||\cdot||_{W^r_p}$.  
\end{definition}
To see more direct definition see \cite{xu} or \cite{BD}. See also \cite[Definition 2,3 and Lemma 2.7]{BDK}.
Despite our characterization also holds for  fractional Sobolev spaces, Theorem \ref{tw:charakteryzacjaBiesova}, we omit these  considerations. 

 \begin{theorem}\label{tw:ORJ}{\color{white}{.}}
 \begin{enumerate}
\item
The operator $E^+$ is localized on $A_+$ and the operator $E^-$ on $A_-$,
\item The both operators $E^\pm$ (i.e., $E^+,E^-$)   $E^\pm:L^2(\mathbb S^d) \to L^2(\mathbb S^d)$  are orthogonal projections

\item For all $r=0,1,\ldots$ and $1\le p<\infty,$ each  $E^\pm$ is a continuous operator 
\[
E^\pm: W^r_p(\Sd) \to  W^r_p(\Sd).
\] 
\end{enumerate}
 \end{theorem}

Consider Daubechies multivariate wavelets. For  $N\geq 2$, let ${}_N\phi$ be a univariate, compactly supported scaling function with support $\supp {}_N\phi=[0,2N-1]$
associated with the compactly supported, orthogonal univariate Daubechies wavelet ${}_N\psi$, see \cite[Section 6.4]{Dau}. Moreover it is known that a smoothness of Daubechies wavelets $\varrho=\varrho(N)\approx 0.2 N$, ${}_N\phi\in C^{\varrho(N)}(\R)$, $\varrho(N)\in \N$,  see \cite[Section 7.1.2]{Dau}. 
  Note that we take ${}_N\psi$ such that $\supp {}_N\psi=[0,2N-1]$.
 For convenience, let $\psi^0={}_N\phi$ and  $\psi^1={}_N\psi$. Let $\E'=\{0,1\}^d$ be the vertices of the unit cube and let $\E=\E' \setminus \{0\} $ be the set of nonzero vertices. For each $\mathbf e=(e_1,\ldots,e_d)\in \E'$, define
\[
\psi^{\mathbf e}(x)=\psi^{e_1}(x_1)\cdots \psi^{e_d}(x_d), \quad x=(x_1,\ldots,x_d)\in \R^d.
\]
Observe that $\supp \psi^{\mathbf e} = [0,2N-1]^d$.

Let $\calD$ be the set of dyadic cubes in $\R^d$ of the form $I=2^{-j}(k+[0,1]^d)$, $j\in \Z$, $k\in \Z^d$. Denote the side length of $I$ by $\ell(I)=2^{-j}$. For any $\mathbf e\in \E'$ define scaled wavelet, related to $I$ by
\[
\psi^{\mathbf e}_I(x)=2^{j d/2}\psi^{\mathbf e}(2^jx-k), \qquad x\in \R^d.
\]
It is well-known that $\{\psi^{\mathbf e}_I(x): I \in \calD, \mathbf e \in \E\}$ is an orthonormal basis of $L^2(\R^d)$, \cite[chapter 3]{Meyer}.  Then to characterize 
classical  Sobolev spaces (as well as  fractional Sobolev
spaces called also Bessel potential spaces)  $W^s_{p}(\R^d)$ for $1<p< \infty$,  $s\geq 0$ or
  Besov spaces $B^s_{pq}(\R^d)$ for $1\leq p\leq\infty$, $1\leq q\leq \infty$, $s>0$
by the magnitude of coefficients of $d$-dimensional Daubechies wavelets  $\psi^e_I$
   we need to assume that
\[
\varrho(N)> s
\]
 \cite[chapter 6, section 2 and section 10]{Meyer}, compare also
 \cite[Theorem 9.4 and Theorem 9.5]{Ker}.
 
  For our purposes it is convenient to consider a localized wavelet systems on a cube.

\begin{figure}[ht]
\centerline{\scalebox{1.1}
       {\includegraphics[scale=0.21]{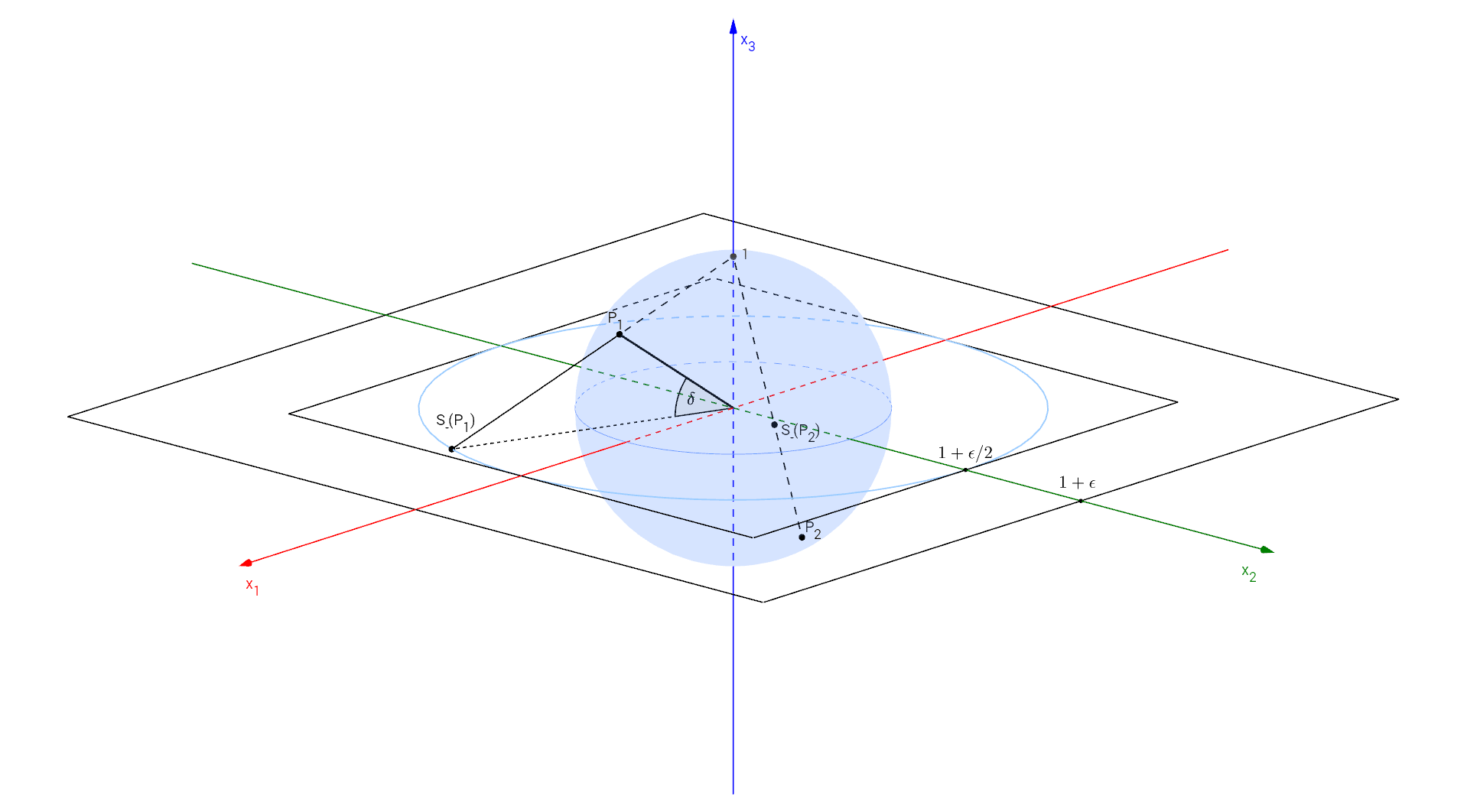}}}
       \caption{Stereographic projection  $S_-$ for $d=2$}
       \label{rys:stereo}
\end{figure}

\begin{definition}\label{def:epsilon}
Suppose that  $J=[-1,1]^d$ is a cube in $\R^d$ and $\ve>0$. Define its $\ve$ enlargement by $J_\ve=[-1-\ve,1+\ve]^d$. Let $j_0 \in \Z$ be the smallest integer such that 
\begin{equation}\label{poziom}
(2N-1)2^{-j_0}\leq  \ve/2.
\end{equation}
For any $j\ge j_0$, consider families of dyadic cubes
\[
\calD_{j}=\calD_{j}(e)=\{ I \in \calD: \ell(I)= 2^{-j} \text{ and } \supp \psi^{\mathbf e}_I \subset J_\ve\} \]
and
\[
\calD_{j_0}^+=\calD_{j_0}^+(e) = \bigcup_{j=j_0}^\infty \calD_j.
\]
Define a localized wavelet system, related to the cube $J$ and  $\ve>0$  by
\begin{equation}\label{los}
S(J,\ve):=\{\psi^{\mathbf e}_{I}:
\mathbf e\in \E, I\in\calD_{j_0}^+(e)\}\cup 
\{\psi^0_{I}:  I \in \calD_{j_0}(0)\}.
\end{equation}
\end{definition}
By appropriate choice of $\ve$, ($\ve=k$  or $\ve=2^{-k}$, $k\in \N\setminus \{0\}$) we get a sequence of finite dimensional spaces 
\begin{equation}\label{MRA}
V_{j_0}\subset \cdots \subset V_j\subset \cdots \subset L^2(J_\ve),
\end{equation}
where
\begin{equation}\label{MRA0}
V_{j}=\spa_{L^2(J_\ve)} \{\psi^0_I : I \in \calD_{j}(0)\}.
\end{equation}
We have usual a dilation and translation properties. We consider only the dilation by two, for $j\geq j_0$  if 
\begin{equation}\label{MRAdilation}
f\in V_j \Rightarrow  f(2\cdot  )\in V_{j+1}.
\end{equation}
We will transport that sequence by two stereographic projections on sphere. After using AWW operators we obtain MRA on $L^2(\Sd)$, (see below).
\begin{lemma}\label{kostka}
The localized wavelet system $S(J,\ve)$ has following properties:
\begin{itemize}
\item
$S(J,\ve)$ is an orthonormal sequence in $L^2=L^2(J_\ve)$,
\item
for every $f\in L^2(J_\ve)$ with $\supp f\subset J_{\ve/2}$ we have 
\begin{equation}\label{framka}
\|f\|^2_{L^2} =\sum_{\mathbf e\in \E} \sum_{I\in \calD_{j_0}^+} |\lan f,\psi^{\mathbf e}_{I}\ran_{L^2}|^2+
\sum_{I\in\calD_{j_0}} |\lan f,\psi^0_{I}\ran_{L^2}|^2.
\end{equation}
\item magnitudes of coefficients $\{|\lan f, g\ran|\}_{g\in S(J,\ve)}$ characterize functions $f \in \mathcal F(\R^d)$
satisfying $\supp f\subset J_{\ve/2}$, where $\mathcal F$  is either the Sobolev space $W^s_p(\R^d)$, $0\leq s<\varrho(N)$, $1<p<\infty$ or the Besov space $B^{s}_{p,q}(\R^d)$, $0<s<\varrho(N)$, $1\leq p,q \leq \infty$.
\end{itemize}
{\emph{Proof is similar to the proof of Lemma 6.1 from \cite{BD} since there is no  restriction  for $p,q=\infty$.}}
\end{lemma}
Localized wavelet system  $S(J,\epsilon)$ is transformed to $\Sd$ by stereographic  projections (see Figure \ref{rys:stereo})
 \begin{align*}
 S_-: \Sd\setminus\{\textbf{1}^d\} \to \R^d,\qquad S_-(x_1,\ldots,x_{d+1}) = \left(\frac{x_1}{1-x_{d+1}},\ldots, \frac{x_{d}}{1-x_{d+1}} \right),
 \end{align*}
  \begin{align*}
 S_+: \Sd\setminus\{-\textbf{1}^d\} \to \R^d,\qquad S_+(x_1,\ldots,x_{d+1}) = \left(\frac{x_1}{1+x_{d+1}},\ldots, \frac{x_{d}}{1+x_{d+1}} \right).
 \end{align*}
 
For  $\epsilon>0$ we define variable change operators (for $+$ and $-$)
\[
T_d^\pm: L^2([-1-\ve,1+\ve]^d) \to L^2(\mathbb S^{d})
\]
given by
\[
T_d^\pm(\psi)(u) =  \frac{\psi( S_\pm (u))}{\sqrt{J_{d}( S_\pm(u))}},\qquad u\in \Sd,
\]
where $J_{d}$ the Jacobian of $S_\pm^{-1}$ 
$$
J_d(x_1,\ldots,x_{d})=\left(\frac{2}{1+x_1^2+\ldots+x_{d}^2}\right)^{d}.
$$
Both operators $T_d^\pm$ are isometric isomorphisms.
This leads us to two a local wavelet system on $\Sd$. Namely,
\begin{align*}
\label{jet2}
{\mathcal S}_\pm&={\mathcal S}_\pm(A_\pm)= T_d^\pm(S(J,\epsilon))= \\
&=\{ T_d^\pm(\psi^{\mathbf e}_{I}):
\mathbf e\in \mathcal{E}, I\in\calD_{j_0}^+\}\cup 
\{ T_d^\pm(\psi^{\mathbf 0}_{I}):  I \in \calD_{j_0}\},
\end{align*}

Various equivalent norms  of Besov spaces $B^{s}_{p,q}=B^{s}_{p,q}(\Sd)$ are given in \cite{LR}. Let us recall 
the definition of Besov space (called Nikolskij-Besov space) form \cite{LR}. For $r\in \N$ let $\omega_r(f,\tau)_p$ be a modulus of smoothness on sphere, i.e.,
\[
\omega_r(f,\tau)_p=\sup_{0<t\leq\tau} \|\Delta^r_t f\|_p
\]
and
\[
\Delta^r_t=(Id-T_t)^r,\quad 0<t<\pi.
\]
Here $Id$ is an identity operator and $T_t$ is a translation operator, compare Definition 2.1.4 \cite{xu},
\[
T_t(f)(\theta)=\frac{\Gamma(d/2)}{2\pi^{d/2}(\sin t)^{d-1}} \int_{\theta\circ y=\cos t}f(y) dl_{\theta,t}(y)
\] 
and $dl_{\theta,t}$ denotes Lebesgue measure on the set $\{y\in \Sd: \theta\circ y=\cos t\}$.
 Let $s>0$ and $1\leq p, q\leq \infty$. Then
\[
B^s_{p,q}(\Sd)=
\{f\in L^p(\Sd): \left(\int_0^\pi \left(\frac{\omega_r(f,t)_p}{t^s}\right)^q \frac{dt}{t} \right)^{1/q}<\infty \},
\]
where $r>s/2$.
We have an analogue of  \cite[Lemma 6.2]{BD}

\begin{lemma}\label{lem:S_+-}
The system ${\mathcal S}_\pm$ has following properties:
\begin{enumerate}
\item
${\mathcal S}_\pm$ is an orthogonal system in  $L^2(\Sd)$,
\item 
If $\epsilon\geq  2\left(\frac{\cos\delta}{1-\sin\delta}-1\right)$, then
$E^\pm(\mathcal S_\pm)$ is   Parseval frame for $E^\pm(L^2(\Sd))$, a.e. for all $f\in E^\pm(L^2(\Sd))$
\begin{equation*}\label{framka2}
\|f\|^2_{2} =\sum_{g\in\mathcal S_\pm} |\lan f,E^\pm(g) \ran|^2.
\end{equation*}

\item the coefficients $\{|\lan f,g \ran|\}_{g\in E^\pm(\mathcal S_\pm)}$ characterize  functions  $f\in E^\pm(\mathcal F(\Sd))$, where  $\mathcal F$ is a Sobolev space $W^s_p(\Sd)$, where $0\leq s<\varrho(N)$, $1<p<\infty$  or Besov space  $B^{s}_{p,q}(\Sd)$, $0<s<\varrho(N)$, $1\leq  p,q \leq\infty$. 
\end{enumerate}
{\emph{Proof is a consequence of Lemma \ref{kostka} and Lemma 6.2 from \cite{BD}.}}
\end{lemma}
We define a wavelet  system called stereographic wavelets corresponding to  $A_+$ and $A_-$
\begin{align}\label{def:ukladS}
\mathcal S:= E^+(\mathcal S_+) \cup E^-(\mathcal S_-),
\end{align}
where 
\[
E^\pm(\mathcal S_\pm)= \{ E^\pm \circ T_d^\pm(\psi^{\mathbf e}_{I}):
\mathbf e\in \mathcal{E}, I\in\calD_{j_0}^+\}\cup 
\{E^\pm\circ T_d^\pm(\psi^{\mathbf 0}_{I}):  I \in \calD_{j_0}\}.
\]
We can define a sequence of finite dimensional spaces $j\geq j_0$
\begin{equation}\label{MRA3}
V_j(\Sd)=\spa_{L^2(\Sd)}  \{E^+\circ T_d^+(\psi^{\mathbf 0}_{I}):  I \in \calD_{j}(0)\}
\oplus \spa_{L^2(\Sd)}\{E^-\circ T_d^-(\psi^{\mathbf 0}_{I}):  I \in \calD_{j}(0)\},
\end{equation}
\[
W_j(\Sd)=\spa_{L^2(\Sd)} \left\{E^\pm\circ T_d^\pm(\psi^{e}_{I}):  e\in \E, I \in \calD_{j} \right\}.
\]
Note that
\[
V_j(\Sd)=E^+\circ T_d^+(V_j)\oplus E^-\circ T_d^-(V_j).
\]
and
\[
W_j(\Sd)=E^+\circ T_d^+(W_j)\oplus E^-\circ T_d^-(W_j).
\]
For a function $f\in V_j(\Sd)$ we will use a notation $f=E^+(f^+)\oplus E^-(f^-)$ where
$f^+ \in  T_d^+(V_j)$ and $f^-\in   T_d^-(V_j).$

\odpar
We have an analogue of Theorem 6.2 from \cite{BD} with multiresolution structure on $L^2(\Sd)$. To formulate a dilation property we need two spherical dilations $\vartheta^\pm$. We take natural parametrization. Let $\Phi_d(\theta,\xi)=x\in \Sd \setminus\{-\textbf{1}^d\}$. Then
\[
\vartheta^+:\Sd\setminus\{-\textbf{1}^d\} \to \Sd\setminus\{-\textbf{1}^d\}
\]
is given by
\[
\vartheta^+(x)=\vartheta^+(\Phi_d(\theta,\xi))=\Phi_d(\phi(\theta),\xi)=(\xi\sin(\phi(\theta)),\cos(\phi(\theta))),
\]
where $\phi:(0,\pi) \to (0,\pi)$ is a diffeomorphism. To find a formula for $\phi=\phi(\theta)$ we need to solve
\[
\frac{2\sin\theta}{1+\cos\theta}=\frac{\sin\phi}{1+\cos\phi}.
\]
Since in the natural parametrization 
\begin{equation}\label{klucze1}
S_+(\xi\sin \theta,\cos \theta)=\frac{\xi\sin\theta}{1+\cos \theta},
\end{equation}
then for the dilation operator  $dial(y)=2y$, where  $y\in \R^d$ we have
\begin{equation}\label{klucze}
S_+ \circ \vartheta^+ =dial\circ  S_+.
\end{equation}
We have a diagram

\begin{center}
\begin{tikzpicture}
  \matrix (m) [matrix of math nodes,row sep=2cm,column sep=4cm,minimum width=2cm]
  {
  \R^d & \R^d \\
   \Sd\setminus\{-\textbf{1}^d\} & \Sd\setminus\{-\textbf{1}^d\} \\
   };
  \path[-stealth]
    (m-2-1) edge node [left] {$S_+$}
 (m-1-1)   (m-1-1)  edge  node [above] {$dial$} (m-1-2)
    (m-2-1.east|-m-2-2) edge node [below] {$\vartheta^+$} (m-2-2)
               (m-2-2) edge node [right] {$S_+$} (m-1-2);
           \end{tikzpicture}
\end{center}

Similar for $\vartheta^-$. 
\begin{theorem}\label{tw:charakteryzacjaBiesova} 
Let $\epsilon=k$ or $\epsilon=2^{-k}$, $k\in \N\setminus \{0\}$. If $\epsilon\geq  2\left(\frac{\cos\delta}{1-\sin\delta}-1\right)$, the
wavelet system  $\mathcal S$ is a Parseval frame in $L^2(\Sd)$. The sequence of $\{V_j(\Sd)\}_{j\geq j_0}$ has following properties
\begin{equation}\label{MRA1}
V_{j_0}(\Sd)\subset V_{j_0+1}(\Sd) \subset \cdots \subset L^2(\Sd),
\end{equation}
\begin{equation}\label{MRA4}
V_{j+1}(\Sd)=V_j(\Sd)\oplus W_j(\Sd),
\end{equation}
and the set
$\bigcup_{j\geq j_0} V_j(\Sd)$ is dense in $L^2(\Sd)$. We have also spherical dilations property, if $j\geq j_0$ and  $f=E^+(f^+) \oplus E^-(f^-)\in V_{j}(\Sd)$, then there are  functions $H^\pm$ independent of $j$ such that
\begin{equation}\label{MRA11}
E^{\pm}\left(H^\pm(\cdot)  f^\pm(\vartheta^\pm(\cdot))\right)\in  V_{j+1}(\Sd).
\end{equation}
Moreover, the magnitudes  of the coefficients $\{|\lan f,g \ran|\}_{g\in \mathcal S}$ characterize $f\in \mathcal F(\Sd)$, where  $\mathcal F$ is Sobolev space $W^s_p(\Sd)$, where $0\leq s<\varrho(N)$, $1<p<\infty$ or Besov space  $B^{s}_{p,q}(\Sd)$, $0<s<\varrho(N)$, $1\leq p,q \leq\infty$. \\
{\emph{[proof in the appendix \ref{A4}]}}
\end{theorem}
In the next section we will consider $\Biess$ for  $0<s<\varrho(N).$ We will use a characterization of functions from $\Biess$ by frame coefficients. Namely, from Theorem \ref{tw:charakteryzacjaBiesova} we get that
a function  $f\in L^2(\Sd)$ belongs to $\Biess,\; 0<s<\varrho(N)$, if $
\|f \|_{s,2}<\infty
$, where
\begin{align}\label{eqdef:charakteryzacja_biesowa}
\|f \|_{s,2}:= \max&\left(  
\sup_{j\geq j_0} 2^{js} \sqrt{\sum_{e\in \mathcal{E}, I\in\calD_{j}}
 \langle f,   E^+ \circ T_d^+(\psi^{\mathbf e}_{I})  \rangle^2},
\; 2^{j_0s} \sqrt{ \sum_{I\in\calD_{j_0}}
 \left\langle f,  E^+ \circ T_d^+(\psi^{\mathbf 0}_{I})  \right\rangle^2},\right.\no 
&\quad\left. 
\sup_{j\geq j_0} 2^{js} \sqrt{\sum_{e\in \mathcal{E}, I\in\calD_{j}}
 \langle f,   E^- \circ T_d^-(\psi^{\mathbf e}_{I})  \rangle^2},
 \; 2^{j_0s} \sqrt{ \sum_{I\in\calD_{j_0}}
 \left\langle f,   E^- \circ T_d^-(\psi^{\mathbf 0}_{I})  \right\rangle^2}
  \right).
\end{align}
In fact we will use an equivalent with the same notation
\begin{align}\label{eqdef:charakteryzacja_biesowa1}
\|f \|_{s,2}:=\max   & \left(
\; 2^{j_0s} \sqrt{ \sum_{I\in\calD_{j_0}}
 \left\langle f,  E^+ \circ T_d^+(\psi^{\mathbf 0}_{I})  \right\rangle^2+
   \sum_{I\in\calD_{j_0}}
 \left\langle f,   E^- \circ T_d^-(\psi^{\mathbf 0}_{I})  \right\rangle^2},\right.\no 
&\left.\quad\sup_{j\geq j_0}  2^{js} \sqrt{\sum_{e\in \mathcal{E}, I\in\calD_{j}}
 \langle f,   E^+ \circ T_d^+(\psi^{\mathbf e}_{I})  \rangle^2
 +\sum_{e\in \mathcal{E}, I\in\calD_{j}}
 \langle f,   E^- \circ T_d^-(\psi^{\mathbf e}_{I})  \rangle^2}\right),
 \;
\end{align}
 We define a family of operators
\begin{align*}
K_j :L^2(\Sd) \to V_j,\quad
K_j(f)(x)=\int_{\Sd} K_j(x,y) f(y) dy,
\end{align*}
where
\begin{align}\label{eqdef:jadro_cale}
&K_j(x,y)=\sum_{I\in  \calD_j} E^\pm\circ T_d^\pm(\psi^{\mathbf 0}_{I})(x) \cdot E^\pm\circ T_d^\pm(\psi^{\mathbf 0}_{I})(y).
\end{align}
If $\{x_n\}, n\in \N$ is the Parseval frame in separable Hilbert space $H$, then
for all $x\in H$ and $I\subset \N$
\begin{equation}\label{frame1}
\|x-\sum_{n\in \N\setminus I} \langle x,x_n \rangle x_n\|^2
=
\|\sum_{n\in  I} \langle x,x_n \rangle x_n\|^2
\leq \sum_{n\in I} \langle x,x_n \rangle^2.
\end{equation}
For completeness of arguments let us prove this inequality. We use the theorem that for Parseval frame  there is orthogonal basis $\{f_n\}$ in $H_1$ such that
$H\subset H_1$ and $x_n=P(f_n)$, where $P$ is an orthogonal projection
\[
P: H_1\to H.
\]
Hence for all $f\in H_1$ and all $I\subset \N$
\[
\|f-\sum_{n \in \N\setminus I}\langle f,f_n \rangle f_n\|^2
=
\|\sum_{n \in  I}\langle f,f_n \rangle f_n\|^2
=\sum_{n\in  I} \langle f,f_n \rangle^2.
\]
If we use orthogonal projection, we obtain
\[
\|P(f-\sum_{n\in \N\setminus I} \langle f,f_n \rangle f_n)\|^2=
\|P(f)-\sum_{n\in \N\setminus I} \langle f,f_n \rangle x_n\|^2
\leq 
\|f-\sum_{n\in \N\setminus I} \langle f,f_n \rangle f_n\|^2.
\]
If we use this for $f=Px=x$, we get 
\[
\|x-\sum_{n\in \N\setminus I} \langle x,Pf_n \rangle x_n)\|^2=
\|x-\sum_{n\in \N\setminus I} \langle x,x_n \rangle x_n\|^2
\leq\sum_{n\in I} \langle x,x_n \rangle^2,
\]
which proves \eqref{frame1}.

Using \eqref{frame1}  we obtain corollary
\begin{corollary}\label{szacowanie}
For all functions from Besov space  $B^{s}_{2,\infty}$, $0<s<\varrho(N)$
and all $j\geq j_0$
\[
\|f-K_jf\|_2\leq 2^{-js} \|f\|_{s,2}.
\]
{\emph{[proof in the appendix \ref{A5}]}}
\end{corollary}


\section{Adaptive estimator of density function}
\label{sec:EstymatorAdaptacyjny}



In this section we present a technical version of Theorem \ref{twierdzenie0}. Note that
\begin{lemma}\label{lem:zalozenie}
For all $N\geq 2$ there is $\D$ such that
\begin{align}
 \forall_{y\in \Sd} \forall_{j\geq j_0} \quad \ints K_j^2(x,y)\dS{x}\leq \D 2^{jd}.
\end{align}
{\emph{[proof in the appendix \ref{A6}]}}
\end{lemma}
One can see that the kernels $K_j(\cdot,\cdot)$ fulfill the conditions from Theorem \ref{tw:TalagrandMOJ}. 
For $f_n(j)(x)=\frac{1}{n}\sum_{i=1}^n K_j(x,X_i)$,
where $K_j(\cdot,\cdot)$ is given by \eqref{eqdef:jadro_cale},
we formulate the analogue of  \cite[Theorem~2]{BN}  (our proof is more precise and gives all needed arguments). The idea of choosing the resolution level is taken from Lepski \cite{lep}.

\begin{theorem}\label{tw:adaptive}
Let $d/2< r<R$ and let   $X_1,\ldots X_n$ be i.i.d.  with density function $f\in \Biess$, where ${r\leq s\leq R}$. We assume that Daubechies
wavelet is smooth enough, i.e., $R<\varrho(N)$.

Let $j_{min}$ and $j_{max}$ be such that  $j_0\leq j_{min}\leq j_{max}$ and
\begin{align}\label{eqdef:jmaxjmin}
j_{min}&=\left\lfloor \frac{\log_2n}{2R+d} \right\rfloor,\qquad
j_{max}=\left\lceil \frac{\log_2n}{2r+d} \right\rceil.
\end{align}

Define $\calJ:=\calJ_n=[j_{min},j_{max}]\cap \N$
and
\begin{align}\label{def:jnhat}
\jn=\min\left\{ j\in\calJ:\forall_{l\in\calJ, l>j} \;\ll{ f_n(j)-f_n(l) }^2 \leq \CS (U \vee 1) \frac{2^{ld}}{n} \right\}, 
\end{align}
where $\CS$ is a constant such that
\begin{align}\label{eq:C(S) twMAIN}
\sqrt{\CS}\geq 2+\M 2\sqrt{\D},
\end{align}
where the constant $\D$is from Lemma  {\rm\ref{lem:zalozenie}} and $\M$ is the constant from theorem {\rm\ref{tw:TalagrandMOJ}} depending on $C_E>0$ and $\D$ for kernels \eqref{eqdef:jadro_cale}.

Then for any $r, R, U>0$ there is $c=c(r,R,U)$ such that for all $s, n$ and $\tilde{B}>1$ if $j_n$ is defined by \eqref{def:jnhat} 
we have 
\begin{align}
\label{eq:SzacowanieDlaAdaptacyjnego}
\sup_{f\in \Sigma(s,\tilde{B}),\|f\|_\infty\leq U} \mathbb{E}\|f_n(j_n)-f\|_2^2\leq c \tilde{B}^{2d/(2s+d)}n^{-2s/(2s+d)}. 
\end{align}
{\emph{[proof in the appendix \ref{A7}]}}
\end{theorem}
In practice constant $C(\mathcal{S})$, from the above Theorem, can be chosen as $$C(\mathcal{S})=\left(2+4\max\{10\sqrt{D_N}/3\ ,\ 4/3+\sqrt{2}\}\sqrt{D_N}\right)^2,$$ where $$D_N=\sup\limits_{j\in\{j_{min},...,j_{max}\}}\sup\limits_{y\in \mathbb{S}^d} \ints K_j^2(x,y)\dS{x}/ 2^{jd}.$$ It can be calculated numerically for $y$ from some grid on $\mathbb{S}^d$. It is also known that the Daubechies wavelets smoothness increases with $N$ approximately like $0.2 N$ (see \cite{Dau} chapter 7). It means that for the estimation one should take $N \approx 5 R$.   
\\

In the proof of the above theorem the following lemma was used.
\begin{lemma}\label{krotko}
Under the above construction we have
\[
\mathbb{E}\ll{f_n(j) - \mathbb{E}f_n(j)}^4\leq 4\left(32\D^2 \sigma^4(j,n) + (\sqrt{2}\cdot 2^{-js} \nb{f})^4\right).
\]
{\emph{[proof in the appendix \ref{A8}]}}
\end{lemma}


\section{Numerical results}\label{sec5}
The first purpose of this section is to explain how one can implement our estimator in $\mathbb{S}^2$, by presenting the exact formula of the estimator with a special choice of functions and parameters. The second purpose is to show, in a numerical experiment, that the estimator works and can be used in practice. 

\begin{figure}[h!]
       {\includegraphics[trim = 30mm 20mm 20mm 20mm,clip, scale=0.88]{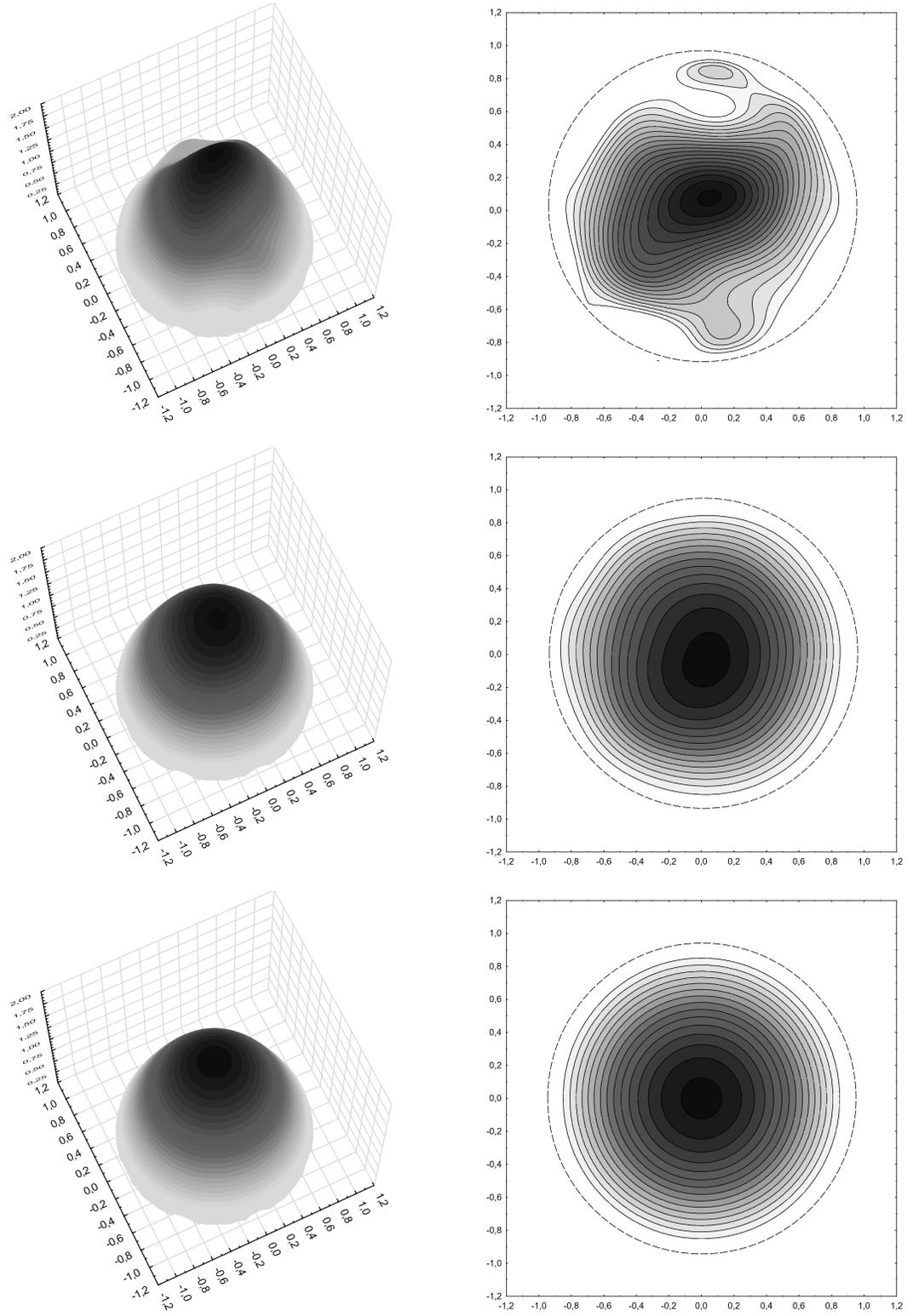}}
       \caption{From top: estimator of function $f_1$ for $n=100$, $n=10000$ and true function $f_1$ on the bottom.}\label{rys:dzwon}
       
\end{figure}

\begin{figure}[h!]
       {\includegraphics[trim = 30mm 20mm 20mm 20mm,clip, scale=0.88]{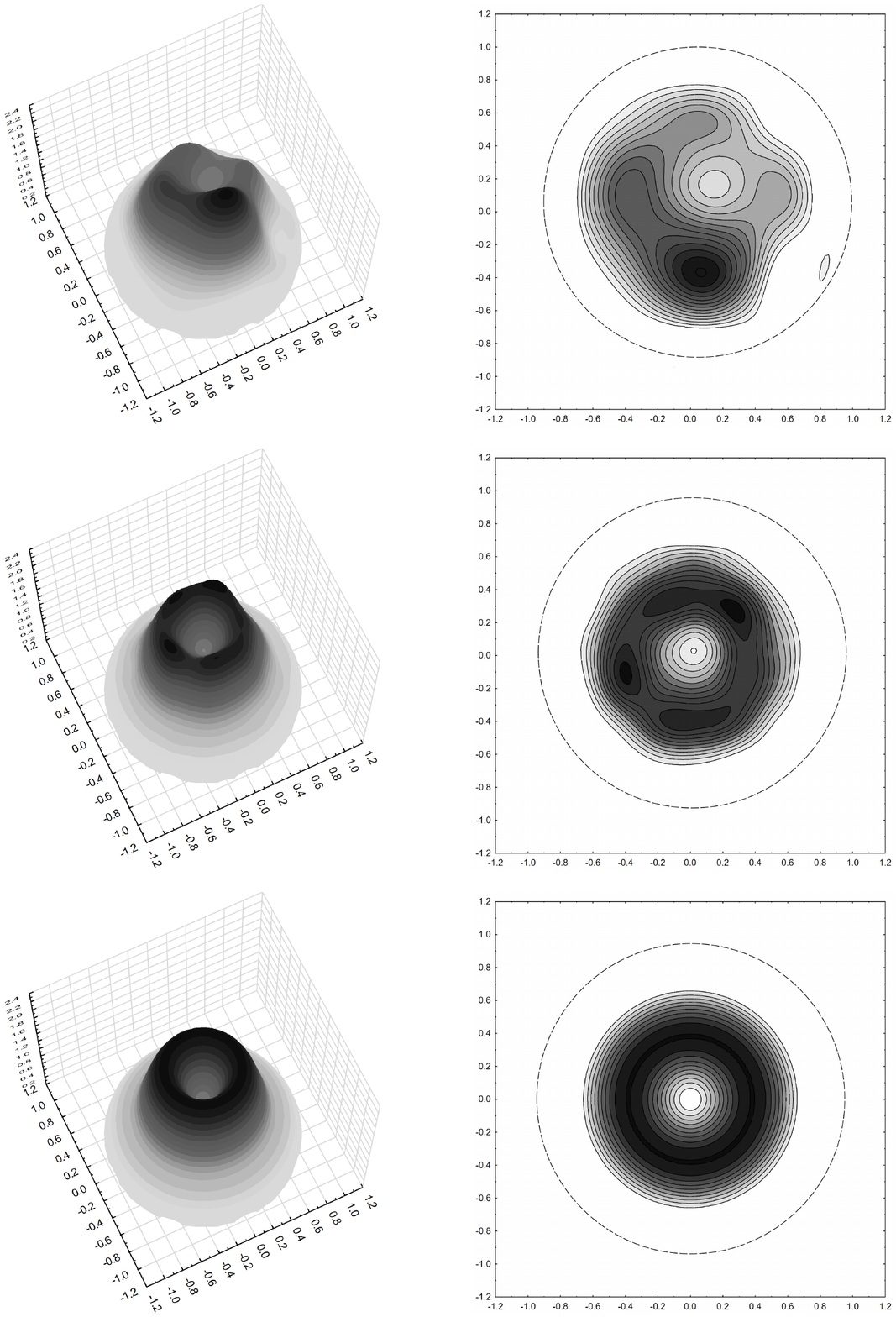}}
       \caption{From top: estimator of function $f_2$ for $n=100$, $n=10000$ and true function $f_2$ on the bottom.}\label{rys:wulkan}
       
\end{figure}

For an estimation the Daubechies wavelets "DB8" with the support $[0,15]$ are used. Two values of
the experiment size are used: $n = 100$ and $n = 10000$. The maximum resolution levels $j_{max}$ from our main Theorem are $j_{max}=2$ for $n=100$ and $j_{max}=3$ for $n=10000$, where $r=3/2$ (see Theorem \ref{tw:adaptive}). Since our wavelet support length is $15$ on the resolution level $0$, then we decide to use minimal resolution level $j_{min}\geq 2$ which allows us to take a lower value for $\epsilon$. Because of that, we set $j_{min}=\max\{2,\left\lfloor \frac{\log_2n}{2R+d} \right\rfloor\}$ which gives us $j_{min}=2$ for $n=100$ and $j_{min}=2$ for $n=10000$, where $R=2$. The resolution levels of the estimator for generated data was $j_n=2$ for $n=100$ and $j_n=3$ for $n=10000$ (the method form the Theorem \ref{tw:adaptive} for choosing $j_n$ can be more useful in practice for bigger sample sizes $n$, when there is a bigger set of possible resolution levels).     
In the estimator formula (see definition \ref{def:estymator}) we choose $\delta=\pi/6$ (see figure \ref{rys:laty}), $\epsilon=4$ (see definition \ref{def:epsilon}) and the distribution function $s \in C^{\infty}(\mathbb{R})$ (see AWW operator (\ref{def:AWW}))
$$s(t)=\left[\ \exp\left(\frac{t-\delta}{t+\delta}\right)\ \Big{/}\ \sqrt{\exp\left(2\ \frac{t-\delta}{t+\delta}\right)+\exp\left(2\ \frac{-t-\delta}{-t+\delta}\right)}\ \right] \mathbb{1}_{(-\delta,\delta)}(t) + \mathbb{1}_{[\delta,\infty)}(t).$$
Notice that the choice of $\epsilon$ is justified since the length of the effective support of DB8 scaling function on the resolution level $2$ is smaller than $\epsilon/2=2$. For the calculation of the wavelets values, a dyadic discretization is used. The distance between discretization points on the resolution level $j$ is $2^{-(j+10)}$.
Data samples $(X_1, Y_1, Z_1),...,(X_n, Y_n, Z_n)$ are generated from the following density functions on the sphere $x^2+y^2+z^2=1$:
\begin{itemize}
\item $f_1(x,y,z)=0,3785\ (\arcsin z - \pi/8)^2\ (\arcsin z - 7\pi/8)^2\ \mathbb{1}_{[\sin(\pi/8),1]}(z)$,
\item $f_2(x,y,z)=42.2126\ (\arcsin z - \pi/4)^2\ (\arcsin z - \pi/2)^2\ \mathbb{1}_{[\sin(\pi/4),1]}(z)$,
\end{itemize}
using the elimination method, which is the following. First we generate an observation $(U_x,U_y,U_z)=(N_x/\sqrt{N^2_x+N^2_y+N^2_z},N_y/\sqrt{N^2_x+N^2_y+N^2_z},N_z/\sqrt{N^2_x+N^2_y+N^2_z})$ from the uniform distribution on the unit sphere, where $N_x, N_y, N_z$ are independent and have standard normal distribution. Then we generate $M$ which is independent of $(U_x,U_y,U_z)$ and has uniform distribution on $[0; 1 + \sup f]$, where $f$ is a density on sphere. If $M<f(U_x,U_y,U_z)$ then we keep the observation $(U_x,U_y,U_z)$. If not, we repeat the procedure. We repeat this until we have the whole sample of size n. The sample is i.i.d. with the density $f$.    
The estimator can be calculated in any point $(x,y,z)$ of the unit sphere by the following formula:
$$\hat{f}_n(j_n)(x,y,z)=\frac{1}{n} \sum\limits_{i=1}^{n} K_j((x,y,z),(X_i,Y_i,Z_i)),$$
where
\begin{eqnarray}
K_j((x,y,z),(X_i,Y_i,Z_i))=\sum\limits_{I\in D_j}& &\left[E^+\circ T^+_d(\psi^{\mathbf 0}_I)(x,y,z)\cdot E^+\circ T^+_d(\psi^{\mathbf 0}_I)(X_i,Y_i,Z_i)\right. \nonumber \\
& & +\left.E^-\circ T^-_d(\psi^{\mathbf 0}_I)(x,y,z)\cdot E^-\circ T^-_d(\psi^{\mathbf 0}_I)(X_i,Y_i,Z_i)\right], \nonumber
\end{eqnarray}
$$
T^+_d(\psi^{\mathbf 0}_I)(x,y,z)=\frac{1}{2}\psi^{\mathbf 0}_I\left(\frac{x}{1+z},\frac{y}{1+z}\right) \left(1+\left(\frac{x}{1+z}\right)^2+\left(\frac{y}{1+z}\right)^2\right),
$$
$$
T^-_d(\psi^{\mathbf 0}_I)(x,y,z)=\frac{1}{2}\psi^{\mathbf 0}_I\left(\frac{x}{1-z},\frac{y}{1-z}\right) \left(1+\left(\frac{x}{1-z}\right)^2+\left(\frac{y}{1-z}\right)^2\right),
$$
\begin{equation}\nonumber
E^-\circ T^-_d(\psi^{\mathbf 0}_I)(x,y,z)=\left\{ \begin{array}{l}
T^-_d(\psi^{\mathbf 0}_I)(x,y,z),\ \ \ \ \ \mathrm{if}\ \ z<\sin\delta\\
\\
s^2(-\arcsin z) T^-_d(\psi^{\mathbf 0}_I)(x,y,z) + s(-\arcsin z) s(\arcsin z)\\ \cdot T^-_d(\psi^{\mathbf 0}_I)(x,y,-z),\  \ \ \ \ \ \mathrm{if}\ \ -\sin\delta \leq   z  \leq \sin\delta \\
\\
0,\ \ \ \ \ \mathrm{if}\  z>\sin\delta\\\\
\end{array}\right. ,
\end{equation}
\\
\\
\begin{equation}\nonumber
E^+\circ T^+_d(\psi^{\mathbf 0}_I)(x,y,z)=\left\{ \begin{array}{l}
0, \ \ \ \ \ \mathrm{if}\ \ z<\sin\delta\\
\\
T^+_d(\psi^{\mathbf 0}_I)(x,y,z)- s^2(-\arcsin z) T^+_d(\psi^{\mathbf 0}_I)(x,y,z) - s(-\arcsin z)\\ \cdot s(\arcsin z) T^+_d(\psi^{\mathbf 0}_I)(x,y,-z),\  \ \ \ \ \ \mathrm{if}\ \ -\sin\delta \leq   z  \leq \sin\delta \\
\\
T^+_d(\psi^{\mathbf 0}_I)(x,y,z),\ \ \ \ \ \mathrm{if}\  z>\sin\delta\\\\
\end{array}\right. .
\end{equation}
\\
\\
In our simulations the estimator values are calculated on the following discrete set of points:
\begin{eqnarray}
\{(x,y,z):& &\hspace{-0.6cm} \left[(x,y)\in \{-0.98,-0.96,...,0.98\}^2 \wedge x^2+y^2 \leq 1 \wedge z^2=1-x^2-y^2 \right]  \nonumber \\ 
& &\hspace{-0.9cm}\vee \left[(x,z)\in \{-0.98,-0.96,...,0.98\}^2 \wedge x^2+z^2 \leq 1 \wedge y^2=1-x^2-z^2 \right] \nonumber \\
& &\hspace{-0.9cm}\vee \left[(y,z)\in \{-0.98,-0.96,...,0.98\}^2 \wedge y^2+z^2 \leq 1 \wedge x^2=1-y^2-z^2 \right]\},
 \nonumber
 \end{eqnarray}
which is approximately uniformly distributed on the sphere and quite comfortable in implementation. 
The results of our estimation are presented in figures \ref{rys:dzwon} and \ref{rys:wulkan}.

\newpage

\newpage

\newpage

\appendix

\section {Mathematical Proofs}
\subsection{Proof of Lemma \ref{lem:SzacowanieWariancji}}\label{A1} 
First we prove (\ref{lSW1}). 
Let $j\geq j_0$ be fixed. For $x\in\Sd$  define
$$Y_i(x) = K_j(x,X_i) - \mathbb{E}K_j(x,X_i)=K_j(x,X_i) - K_jf(x).$$
Since  $Y_i(x)$ are i.i.d. and $\mathbb{E}Y_i(x) = 0$ by \eqref{war2} we get 
\begin{align*}
\mathbb{E}\ll{f_n(j)-\mathbb{E}f_n(j)}^2
&= \frac{1}{n}\ints \mathbb{E}\left( Y_i(x)\right)^2\dS{x}\nonumber\\
&\leq \frac{1}{n}\ints\ints  K_j^2(x,y)f(y)\dS{y}\dS{x} \leq D\frac{2^{jd}}{n}.
\end{align*}
Now we prove (\ref{lSW2}). From Jensen's inequality and \eqref{lSW1} we have
\begin{align*}
 \left(\mathbb{E}\ll{\sum_{i=1}^n \Big(K_j(\cdot,X_i)- \mathbb{E}K_j(\cdot,X_i)\Big)}\right)^2
  &\leq n^2 \mathbb{E}\ll{\frac{1}{n}\sum_{i=1}^n \Big(K_j(\cdot,X_i)- \mathbb{E}K_j(\cdot,X_i)\Big)}^2\\
  &= n^2 \mathbb{E}\ll{f_n(j)-\mathbb{E}f_n(j)}^2 \leq Dn2^{jd}.
\end{align*}

\subsection{Proof of Lemma \ref{assumption}}\label{AA}
Let us check that $\Cal{K}$ satisfies the assumptions of Theorem \ref{tw:talagrand}. We start with
assumption that $k_g\in \K$ are uniformly bounded. If we take $k_g\in\Cal{K}$, then
\begin{align*}
\|k_g\|_\infty&=\sup_{x\in\Sd}|k_g(x)|\\ 
&\leq \sup_{x\in\Sd} (\ll{g}\ll{K_j(\cdot,x)}+\ll{g}\ll{K_j(f)})\kom{Schwarz's inequality}\\
&\leq \sup_{x\in\Sd} (\ll{K_j(\cdot,x)}+\ll{K_j(f)})\kom{since $g\in B_0$}\\
&\leq \sup_{x\in\Sd} (\ll{K_j(\cdot,x)}+\ll{f})\kom{ (\ref{war3}.)}.\\
\end{align*}
Note that by (\ref{war2}) we get
\begin{eqnarray}\label{raz}
\forall_{j\geq j_0}\forall_{x\in\Sd}\quad\ll{K_j(\cdot,x)}&=& \left( \ints K_j^2(t,x) \dS{t}\right)^{1/2}\leq \sqrt{D}2^{jd/2}.
\end{eqnarray}
Moreover
\begin{eqnarray}\label{dwa}
\ll{f}&=& \left( \ints f^2(t) \dS{t}\right)^{1/2}\leq  \left(\|f\|_\infty \ints f(t) \dS{t}\right)^{1/2} = \|f\|_\infty^{1/2} < \infty.\nonumber\\
\end{eqnarray}
By (\ref{raz}) and (\ref{dwa})
\begin{eqnarray*}
\|k_g\|_\infty&\leq& \sqrt{D}2^{jd/2} + \|f\|_\infty^{1/2} =: U_{\Cal{K}_j}<\infty.
\end{eqnarray*}
We have proved that a function from $\Cal{K}=\Cal{K}_j$ is  bounded by $U_{\Cal{K}_j}.$
By (\ref{eqdef:U_K}) we get  
$$
\forall_{k_g\in \Cal{K}}\quad\ints |k_g(x)|f(x)\dS{x}<\infty.
$$
Consequently,
\begin{eqnarray*}
\mu(k_g) &=& \ints k_g(x)f(x)\dS{x}\nonumber\\ 
&=&\ints g(t) \left(  K_j(f)(t) - K_j(f)(t) \ints f(x)\dS{x} \right) \dS{t} = 0. \nonumber
\end{eqnarray*}
Let us show \eqref{szacowanie} i.e., for   $\omega_{\Cal{K}}^2=\|f\|_\infty$  
$$
\forall_{g\in B_0}\quad \mathbb{E}\big[k_g(X)\big]^2\leq \omega_{\Cal{K}}^2.
$$
If $g\in B_0$, then by \eqref{war3}
\begin{align}
\mathbb{E}\big[k_g(X)\big]^2&= Var\left(\ints g(t) K_j(t,X)\dS{t} \right)\nonumber\\
&\leq \mathbb{E}\left[\ints g(t) K_j(t,X)\dS{t} \right]^2\nonumber\\
&\leq  \|f \|_\infty\ints\left(\ints g(t) K_j(t,x)\dS{t} \right)^2\dS{x}\nonumber\\
&\leq \|f \|_\infty \ll{g}^2.
\end{align}
Consequently,
\begin{eqnarray*}
\sup_{g\in B_0\subset B} \mathbb{E}\big[k_g(X)\big]^2  
&=&\|f \|_\infty \sup_{g\in B_0\subset B}\ll{g}^2\leq \omega_{\Cal{K}}^2.
\end{eqnarray*}

\subsection{Proof of Theorem \ref{tw:charakteryzacjaBiesova}}\label{A4}
From Lemma \ref{lem:S_+-} we get that $\mathcal S$ is a Parseval frame in $L^2(\Sd)$. From \eqref{MRA} we obtain \eqref{MRA1}. Since  $\mathcal S$ is a Parseval frame in $L^2(\Sd)$ we get that the sum $\bigcup_{j\geq j_0} V_j(\Sd)$ is dense in $L^2(\Sd)$. 
Let 
\[
W_j=\spa_{L^2(J_\ve)} \{ \psi^e_I: e\in \E, I\in \calD_j(e) \}.
\]
Then by definition of $W_j$ and \eqref{MRA0}
\[
V_j\oplus W_j \subset V_{j+1}.
\]
Hence by \eqref{MRA3}
\[
V_j(\Sd)+ W_j(\Sd)\subset V_{j+1}(\Sd).
\]
Note that both orthogonal projection $E^+$ and $E^-$ are localized.
 Then if $f\in L^2(J_\ve)$ is such that 
\[
\supp f \cap (-1-\ve/2,1+\ve/2)^d =\emptyset
\]
then
\[
E^\pm\circ T_d^\pm (f)=0.
\]
Hence the space $V_{j+1}(\Sd)$ is spanned by functions $E^\pm\circ T_d^\pm\psi_L^0$ such that
$l(L)=2^{-j-1}$ and 
\[
\supp \psi_L^0 \cap (-1-\ve/2,1+\ve/2)^d \neq \emptyset.
\]
From definition of $j_0$ we get that such $\psi_L^0$ are 
\[
\psi_L^0=\sum_{I\in \calD_j} \langle \psi_L^0, \psi_I^0 \rangle \psi_I^0+
\sum_{e\in \E, I\in \calD_j} \langle \psi_L^0, \psi_I^e \rangle \psi_I^e.
\] 
Consequently
\[
V_{j+1}(\Sd)\subset V_j(\Sd)+ W_j(\Sd)
\]
and we prove \eqref{MRA1}. Now we  prove \eqref{MRA11}. Let us take $\psi_I^0$ and $I\in \calD_j(0)$. Then 
\[
 T_d^+(\psi_I^0)(\vartheta^+ (u))=\frac{\psi_I^0( S_+ (\vartheta^+ (u))}{\sqrt{J_{d}( S_+ (\vartheta^+ (u)))}},\qquad u\in \Sd,
 \]
 where
$$
J_d(x_1,\ldots,x_{d})=\left(\frac{2}{1+x_1^2+\ldots+x_{d}^2}\right)^{d}.
$$
By \eqref{klucze} we have
\[
 T_d^+(\psi_I^0)(\vartheta^+ (u))=\frac{\psi_I^0( 2 S_+  (u))}{\sqrt{J_{d}( 2 S_+  (u))}},\qquad u\in \Sd.
 \]
 Consequently there is $\tilde I\in \calD_{j+1}(0) $ such that for $y\in \R^d$
 \[
 \psi_I^0( 2 y)=\psi_{\tilde I}^0 (y).
 \]
 The function $H^+$ we obtain comparing 
 $ T_d^+(\psi_I^0)(\vartheta^+ (u))$ with  
\[
 T_d^+(\psi_{\tilde I}^0)(u)=\frac{\psi_{\tilde I}^0( S_+ (u))}{\sqrt{J_{d}( S_+ (u))}},\qquad u\in \Sd.
 \]  
 Hence
 \[
 H^+(u)=\sqrt{\frac{J_{d}( 2 S_+ (u))}{J_{d}( S_+ (u))}}.
 \]
By \eqref{klucze1} in the natural parametrization $u=(\xi\sin \theta,\cos\theta)$ the function $H^+$ depends on $a=a(\theta)=\sin\theta/(1+\cos\theta)$ since  
\[
J_d(S_+(u))=J_d\left(\frac{\xi\sin\theta}{1+\cos \theta}\right)=J_d(a\xi)=\left( \frac{1}{1+a^2} \right)^d.
\]  
From Theorem~6.1 ~\cite{BD} we get that there is  $C>0$ such that
for  $f\in \mathcal F(\Sd).$
$$
\| E^\pm f \|_{\mathcal{F}(\Sd)} \leq C \|  f \|_{\mathcal{F}(\Sd)}. 
$$ 
Since both operators  $E^+$ i $E^-$ are projection in Banach spaces, then
$$
\mathcal{F}(\Sd) = E^+ (\mathcal{F}(\Sd))\oplus E^- (\mathcal{F}(\Sd))
$$
with equivalence of norm
$$
\| f \|_{\mathcal{F}(\Sd)} \simeq \|E^+f \|_{\mathcal{F}(\Sd)} + \|E^-f \|_{\mathcal{F}(\Sd)}.
$$
Since $E^+f \in E^+ (\mathcal{F}(\Sd))$ a $E^-f \in E^- (\mathcal{F}(\Sd)),$ consequently from Lemma \ref{lem:S_+-} we get the Theorem.

\subsection{Proof of Corollary \ref{szacowanie}}\label{A5}
Let
\[
K_{j+1}(f)=K_{j+1}^+(f)+K_{j+1}^-(f),
\]
where
\[
K_{j+1}^+(f)=\sum_{I\in \calD_{j+1}}\langle f,E^+T_d^+(\psi_I^0) \rangle E^+T_d^+(\psi_I^0).
\]
Then
\[
K_{j+1}^+(f)=E^+\left(\sum_{I\in \calD_{j+1}}\langle E^+f,T_d^+(\psi_I^0) \rangle T_d^+(\psi_I^0)\right)
\]
But  for $j\geq j_0$ we have 
\[
V_{j}\oplus W_{j}\subset V_{j+1}.
\]
So since 
\[
supp \left((T_d^+)^{-1}E^+f\right) \subset [-1-\ve/2,1+\ve/2]^d
\]
 we get
\[
K_{j+1}^+(f)=E^+\left(\sum_{I\in \calD_j}\langle E^+f,T_d^+(\psi_I^0) \rangle T_d^+(\psi_I^0)\right)
+E^+\left(\sum_{e\in \E} \sum_{I\in \calD_j(e)}\langle E^+f,T_d^+(\psi_I^e) \rangle T_d^+(\psi_I^e)\right)
\]
If we introduce
\[
Q_j(f)=\int_{\Sd} G_j(x,y) f(y) \sigma_d(dy)
\]
where
\begin{equation}\label{MRA5}
G_j(x,y)=
\sum_{e\in \mathcal{E}}
\sum_{I\in  \calD_{j}} E^\pm\circ T_d^\pm(\psi^{e}_{I})(x) \cdot E^\pm\circ T_d^\pm(\psi^{e}_{I})(y).
\end{equation}
then for $j\geq j_0$
\begin{equation}\label{jadra}
K_{j+1}=K_{j}+Q_j.
\end{equation}
Hence
\begin{equation}\label{jadra1}
K_{j+1}=K_{j_0}+\sum_{k=j_0}^{j}Q_k.
\end{equation}
We get the corollary by \eqref{frame1}, \eqref{eqdef:charakteryzacja_biesowa} and Theorem \ref{tw:charakteryzacjaBiesova}.

\subsection{Proof of Lemma \ref{lem:zalozenie}}\label{A6}
Since $E^\pm$ are orthogonal projections and $T_d^\pm$ are isometric isomorphisms we have
\begin{align}\label{eq:1 lem:zalozenie}
\ll{E^\pm\circ T_d^\pm(\psi^{\mathbf 0}_{I})}\leq \ll{T_d^\pm(\psi^{\mathbf 0}_{I})}=\ll{\psi^{\mathbf 0}_{I}} = 1.
\end{align}
Let
\begin{align*}
 \tilde{E}^+_I(x,y) &:=E^+\circ T_d^+(\psi^{\mathbf 0}_{I})(x) \cdot E^+\circ T_d^+(\psi^{\mathbf 0}_{I})(y)
\end{align*}
and
\begin{align*}
 \tilde{E}^-_I(x,y) &:=E^-\circ T_d^-(\psi^{\mathbf 0}_{I})(x) \cdot E^-\circ T_d^-(\psi^{\mathbf 0}_{I})(y).
\end{align*}
Then for $y\in\Sd$ and $j\geq j_0$
\begin{align}\label{eq:2 lem:zalozenie}
\ints K_j^2(x,y)\dS{x}&= \ints \left[  \sum_{I\in  \calD_j} \left(  \tilde{E}^+_I(x,y) +  \tilde{E}^-_I(x,y)\right)  \right]^2 \dS{x}.
\end{align}
Note that there is a constant $\tilde{D}$ such that for all $x,y\in \Sd$
\begin{align*}
D(x,y)&:=\#\{I\in  \calD_j: \tilde{E}^+_I(x,y) +  \tilde{E}^-_I(x,y)\neq 0\}\leq \tilde{D}<\infty.
\end{align*}
By Jensen's inequality \eqref{eq:2 lem:zalozenie} 
\begin{align}\label{eq:pom1}
\ints K_j^2(x,y)\dS{x}&\leq \ints  \tilde{D}\sum_{I\in  \calD_j} \left[  \tilde{E}^+_I(x,y) +  \tilde{E}^-_I(x,y)\right]^2 \dS{x}\no 
&\leq 2 \tilde{D}\ints  \sum_{I\in  \calD_j} \left[ \left(  \tilde{E}^+_I(x,y)\right)^2 +  \left(\tilde{E}^-_I(x,y)\right)^2  \right] \dS{x}.
\end{align}
From \eqref{eq:1 lem:zalozenie} we get
\begin{align}\label{eq:pom2}
\ints  \sum_{I\in  \calD_j}  \left(  \tilde{E}^\pm_I(x,y)\right)^2 \dS{x} &=   
                   \sum_{I\in  \calD_j}  \Big( E^\pm\circ T_d^\pm(\psi^{\mathbf 0}_{I})(y)\Big)^2 \ints \Big( E^\pm\circ T_d^\pm(\psi^{\mathbf 0}_{I})(x)\Big)^2 \dS{x}\no 
&\leq  \sum_{I\in  \calD_j}  \Big( E^\pm\circ T_d^\pm(\psi^{\mathbf 0}_{I})(y)\Big)^2 \leq C_1\left(2^{jd/2}\right)^2 = C_12^{jd}.
\end{align}
Consequently from \eqref{eq:pom1} i \eqref{eq:pom2}
\begin{align*}
\ints K_j^2(x,y)\dS{x}&\leq 4\tilde{D}C_1\cdot 2^{jd} = \D2^{jd}.
\end{align*}

\subsection{Proof of Theorem \ref{tw:adaptive}}\label{adaptive_proof}\label{A7}
From Lemma  \ref{lem:SzacowanieWariancji} we get that for all $j\geq j_0$ 
\begin{align}\label{def:sigma}
\mathbb{E}\ll{f_n(j)-\mathbb{E}f_n(j)}^2\leq \D \frac{2^{jd}}{n}=:\D \sigma^2(j,n). 
\end{align}
If $f\in \Biess$ by \eqref{eqdef:charakteryzacja_biesowa} and Corollary \ref{szacowanie} we get
\begin{align}\label{def:bias}
\ll{\mathbb{E}f_n(j)-f}^2\leq (\sqrt{2}\cdot 2^{-js} \nb{f})^2=:B^2(j,f). 
\end{align}
We define
\begin{align}\label{def:jstar}
j^*=\min\left\{ j\in\calJ: B(j,f)\leq \sqrt{\D}\sigma(j,n) \right\}. 
\end{align}
Note that
\begin{align}\label{eq:DwieCzesci}
\mathbb{E} \ll{f_n(\jn) - f }^2  &= \mathbb{E}\left( \ll{f_n(\jn) - f }^2   \ch_{\{\jn\leq j^*\}}  + \ll{f_n(\jn) - f }^2 \ch_{\{\jn>j^*\}}   \right).
\end{align}
To show  \eqref{eq:SzacowanieDlaAdaptacyjnego} we will estimate both components. 

\odpar
(I) For first component we use $(a+b)^2\leq$ $2(a^2+b^2)$, thus
\label{dla jn^ <= j*}
\begin{align*}
\mathbb{E} \ll{f_n(\jn) - f }^2\ch_{\{\jn\leq j^*\}}   &\leq 
\mathbb{E} \left(\ll{f_n(\jn) - f_n(j^*)}+\ll{f_n(j^*)-f}\right)^2 \ch_{\{\jn\leq j^*\}}  \\
&\leq 
2 \mathbb{E}\ll{f_n(\jn) - f_n(j^*)}^2\ch_{\{\jn\leq j^*\}}+2\mathbb{E}\ll{f_n(j^*)-f}^2 \ch_{\{\jn\leq j^*\}} \numberthis \label{dwaa} 
\end{align*}

If $\jn\leq j^*$, then by \eqref{def:jnhat} we get
\begin{align*}
\mathbb{E} \ll{f_n(\jn) - f_n(j^*) }^2 \ch_{\{\jn\leq j^*\}} 
& \leq \mathbb{E}\left(  \CS (||f||_\infty \vee 1) \frac{2^{j^*d}}{n}\ch_{\{\jn\leq j^*\}}\right)\no 
& = \CS (||f||_\infty \vee 1) \sigma^2(j^*,n) \no 
\end{align*}

Applying  \eqref{def:sigma}, \eqref{def:bias} and \eqref{def:jstar} we get
\begin{align*}
\mathbb{E} \ll{f_n(j^*) - f }^2  \ch_{\{\jn\leq j^*\}} 
&=      \mathbb{E} \ll{f_n(j^*) - \mathbb{E}f_n(j^*) }^2 + \ll{\mathbb{E}f_n(j^*) - f }^2\\
&\leq   \D\sigma^2(j^*,n) + B^2(j^*,f)
\leq 2\D\sigma^2(j^*,n).\numberthis \label{razz}
\end{align*}
\odpar
Note that for $f\in\Sigma(s,\tilde{B}),$ 
\begin{align}\label{eq:z j*-}
\sigma^2(j^*,n)\leq D' \tilde{B}^{\frac{2d}{2s+d}} n^{\frac{-2s}{2s+d}}.
\end{align}
Indeed 
for  $j^{*-} = j^*-1 $ from \eqref{def:jstar} we have
\begin{align*}
B(j^{*-},f) \geq \sqrt{\D} \sigma(j^{*-},n).
\end{align*}
Consequently using \eqref{def:sigma} and \eqref{def:bias} we get
\begin{align*}
 \nb{f}^2 2^{-2sj^{*-}+1} &\geq \D\frac{2^{dj^{*-}}}{n}\\
\nb{f}^2 2^{-2s(j^{*}-1)+1} &\geq \D\frac{2^{d(j^{*}-1)}}{n}\\
\frac{2}{\D} \nb{f}^2 2^{2s+d} &\geq \frac{1}{n}2^{j^*(2s+d)}.
\end{align*}
Hence by standard calculation we get
\begin{align*}
\left(\frac{2}{\D} \nb{f}^{2}\right)^d 2^{d(2s+d)} &\geq \left(\frac{2^{j^*d}}{n}\right)^{2s+d} n^{2s}.
\end{align*}
From  \eqref{def:sigma}
\begin{align*}
\left(\frac{2}{\D} \nb{f}^{2}\right)^{   \frac{d}{2s+d}   } 2^{d} n^{ \frac{-2s}{2s+d} } &\geq \sigma^2(j^*,n).
\end{align*}
Finally for $f\in\Sigma(s,\tilde{B})$ and $s\geq d/2$
\begin{align*}
\sigma^2(j^*,n) &\leq \left(\frac{2}{\D}\right)^{ \frac{d}{2s+d} }2^{d} \cdot \tilde{B}^{ \frac{2d}{2s+d} }  n^{ \frac{-2s}{2s+d} }
\leq \left(\frac{2}{\D}\right)^{1/2}2^{d} \cdot \tilde{B}^{ 2\frac{d}{2s+d} }  n^{ \frac{-2s}{2s+d} }.
\end{align*}
Note that the constant in the inequality \eqref{eq:z j*-} is the following 
\begin{align}\label{eqdef:D_prim}
D' = \sqrt{\frac{2}{\D}} 2^{d}.
\end{align}
By  \eqref{dwaa}, \eqref{razz} and  \eqref{eq:z j*-} for $
C''=\CS (||f||_\infty \vee 1).
$
\begin{align}\label{KoncowaDlaI}
\mathbb{E} \ll{f_n(\jn) - f }^2\ch_{\{\jn\leq j^*\}}
&\leq 2(C''+\D)\sigma^2(j^*,n)
\leq 2(C''+\D)D' \tilde{B}^{\frac{2d}{2s+d}} n^{\frac{-2s}{2s+d}},
\end{align}

(II) For second component of  \eqref{eq:DwieCzesci} using Schwartz inequality we get
\begin{align}\label{poczatekII}
\mathbb{E}\left( \ll{f_n(\jn) - f }^2 \ch_{\{\jn>j^*\}}\right)&=
\sum_{j\in\calJ,j>j^*} \mathbb{E}\left( \ll{f_n(j) - f }^2 \ch_{\{\jn=j\}}\right) \no 
&\leq \sum_{j\in\calJ,j>j^*} \left(\mathbb{E} \ll{f_n(j) - f}^4\right)^{1/2} \sqrt{P(\jn=j)}.
\end{align}
Note that 
\begin{align*}
\ll{f_n(j) - f}
&\leq \ll{f_n(j) - \mathbb{E}f_n(j)} + \ll{ \mathbb{E}f_n(j) - f}
\leq \ll{f_n(j) - \mathbb{E}f_n(j)} + B(j,f).
\end{align*}
Hence 
\begin{align}\label{eq:E()^4 proof:main}
\mathbb{E}\left(\ll{f_n(j) - f}\right)^4
&\leq  \mathbb{E}\left(\ll{f_n(j) - \mathbb{E}f_n(j)} + B(j,f)\right)^4\no 
&\leq 8 \left(\mathbb{E}\ll{f_n(j) - \mathbb{E}f_n(j)}^4 + B^4(j,f)\right).  
\end{align}
By Lemma \ref{krotko} and \eqref{def:jstar} for $j>j^*$ and $j\in\calJ$ there is $C$ such that
\begin{equation}\label{dodaj}
\mathbb{E} \ll{f_n(j) - \mathbb{E}f_n(j)}^4 \leq  4\left(32\D^2 \sigma^4(j,n) + (\sqrt{2}\cdot 2^{-js} \nb{f})^4\right)\leq C \sigma^4(j_{\max},n).
\end{equation}
Consequently by \eqref{poczatekII}, \eqref{dodaj} and the definition of $j_{\max}$ there is $C''$ such that
\begin{align}\label{eq:part2}
\mathbb{E}\left( \ll{f_n(\jn) - f }^2 \ch_{\{\jn>j^*\}}\right)&= \sqrt{C''} \sum_{j\in\calJ:j>j*} \sqrt{P(\jn=j)}.
\end{align}\\
Let us fix  $j\in \calJ$ such that $j>j*$. We will estimate $P(\jn=j).$
Since $\jn = j,$ then from the definition of $\jn$ and for $j^-=j-1$ we get
that there is $l>j-1\geq j*$ such that
\begin{align}\label{eq:A}
 \quad\ll{f_n(j^-) - f_n(l)}^2 &> \CS (||f||_\infty \vee 1) \frac{2^{ld}}{n}.
\end{align}
Note that 
\begin{align}\label{eq:B1}
&\ll{f_n(j^-) - f_n(l)}\leq\no 
& \leq \ll{f_n(j^-) - \mathbb{E}f_n(j^-) - f_n(l)+ \mathbb{E}f_n(l)} + \ll{f- \mathbb{E}f_n(j^-)} +  \ll{f- \mathbb{E}f_n(l)}\nonumber\\
&\leq \ll{f_n(j^-) - \mathbb{E}f_n(j^-) - f_n(l)+ \mathbb{E}f_n(l)} + B(j^-,f)+ B(l,f).
\end{align}
From the definitions of $B(l,f)$ and $\sigma(l,n)$ we get
$$
B(j^-,f)+ B(l,f) \leq 2B(j^*,f)\leq 2\sqrt{\D}\sigma(j^*,f)\leq 2\sqrt{\D}\sigma(l,f).
$$
By \eqref{eq:B1} 
\begin{align}\label{eq:B2}
\ll{f_n(j^-) - f_n(l)}&\leq \ll{f_n(j^-) - \mathbb{E}f_n(j^-) - f_n(l)+ \mathbb{E}f_n(l)} + 2\sqrt{\D}\sigma(l,f).
\end{align}
Hence by  \eqref{eq:A} and \eqref{eq:B2} we have
\begin{align}\label{Pjn=jI}
&P(\jn=j)\leq\no
&\leq  \sum_{l\in \calJ:l\geq j} P\left(\ll{f_n(j^-) - f_n(l)} > \sqrt{\CS (||f||_\infty \vee 1)} \sigma(l,n)\right)\no 
 &\leq  \sum_{l\in \calJ:l\geq j} P\Big(\ll{f_n(j^-) - \mathbb{E}f_n(j^-) - f_n(l)+ \mathbb{E}f_n(l)}> \no 
  &\hspace{6cm}>\Big(\sqrt{\CS (||f||_\infty \vee 1)} -2\sqrt{\D}\Big)\sigma(l,n)\Big)\no 
 &\leq  \sum_{l\in \calJ:l\geq j} \Big[ P\Big(\ll{f_n(j^-) - \mathbb{E}f_n(j^-)} > G\sigma(l,n)\Big)  
 + P\Big(\ll{f_n(l) - \mathbb{E}f_n(l)} > G\sigma(l,n)\Big)  \Big],
\end{align}
where the constant G is given by
$$
G=\frac{\sqrt{\CS (||f||_\infty \vee 1)} -2\sqrt{\D}}{2}.
$$
Since we assume that $\sqrt{\CS}\geq 2\M+2\sqrt{\D}$, see \eqref{eq:C(S) twMAIN} we get
\begin{align*}
G&=\frac{1}{2}\left(\sqrt{\CS (||f||_\infty \vee 1)} -2\sqrt{\D}\right)\\
&\geq \frac{1}{2}\left(2\M\sqrt{ (||f||_\infty \vee 1)} +2\sqrt{\D(||f||_\infty \vee 1)} -2\sqrt{\D}\right) \\
&\geq   \M\sqrt{ (||f||_\infty \vee 1)}\kom{ $(||f||_\infty \vee 1) \geq 1$}
\end{align*}
Thus, it holds
\begin{align*}
P(\jn=j)&\leq   
\sum_{l\in \calJ:l\geq j} \Big[P\left(\ll{f_n(j^-) - \mathbb{E}f_n(j^-)} > \M\sqrt{ (||f||_\infty \vee 1)}\sigma(l,n)\right)  \\
& \quad\qquad+ P\left(\ll{f_n(l) - \mathbb{E}f_n(l)} > \M\sqrt{ (||f||_\infty \vee 1)}\sigma(l,n)\right) \Big].
\end{align*}
Since for  $l\geq j>j^-$ we have $\sigma(l,n)>\sigma(j^-,n)$ then
\begin{align*}
P(\jn=j)&\leq   
\sum_{l\in \calJ:l\geq j} \Big[P\left(\ll{f_n(j^-) - \mathbb{E}f_n(j^-)} > \M\sqrt{ (||f||_\infty \vee 1)}\sigma(j^-,n)\right)  \\
& \quad\qquad+ P\left(\ll{f_n(l) - \mathbb{E}f_n(l)} > \M\sqrt{ (||f||_\infty \vee 1)}\sigma(l,n)\right) \Big].
\end{align*}
By Theorem \ref{tw:TalagrandMOJ} 
\begin{align*}
P(\jn=j)&\leq  \sum_{l\in \calJ:l\geq j} \left( 2 e^{-2^{dj^-}} + 2 e^{-2^{dl}} \right)\leq 4(j_{max} - j_{min})e^{-2^{dj_{min}}}.
\end{align*}
Using this in \eqref{eq:part2} we obtain
\begin{align}\label{eq:wzor_jmax_jmin}
\mathbb{E}\left( \ll{f_n(\jn) - f }^2 \ch_{\{\jn>j^*\}}\right)&=
\sqrt{C''} \sum_{j\in\calJ:j>j*} \sqrt{4(j_{max} - j_{min})}e^{-2^{dj_{min}}/2}\no
&\leq 2\sqrt{C''} (j_{max} - j_{min})^{\frac{3}{2}} e^{-\frac{1}{2} 2^{dj_{min}}}.
\end{align}
But from \eqref{eqdef:jmaxjmin}
$$
j_{max} - j_{min}\leq \left(\frac{1}{2r+d}-\frac{1}{2R+d}\right)\log_2n
\leq \left(\frac{1}{2r+d}-\frac{1}{2R+d}\right)n
$$
and
$$
e^{-\frac{1}{2}2^{dj_{min}}}\simeq e^{-\frac{1}{2}n^{\frac{d}{2R+d}}}.
$$
Now  \eqref{eq:wzor_jmax_jmin} we estimate by
\begin{align}\label{eq:wzor_jmax_jmin2}
\mathbb{E}\left( \ll{f_n(\jn) - f }^2 \ch_{\{\jn>j^*\}}\right)
&\leq
\tilde{C} \left(\frac{1}{2r+d}-\frac{1}{2R+d}\right)^{\frac{3}{2}}
n^{\frac{3}{2}}e^{-\frac{1}{2}n^{\frac{d}{2R+d}}}\no 
&=\tilde{C_1}\frac{n^{\frac32+\frac{2s}{2s+d}}              }{\exp(\frac12 n^{\frac{d}{2R+d}})}\cdot n^{-\frac{2s}{2s+d}}.
\end{align}
But for any $\alpha, \beta, \gamma >0$ we have
$$
\lim_{x\to \infty} \frac{x^{\beta}}{e^{\alpha x^\gamma}}=0.
$$
Consequently  \eqref{eq:wzor_jmax_jmin2} we estimate by
\begin{align}\label{KoniecDlaII}
\mathbb{E}\left( \ll{f_n(\jn) - f }^2 \ch_{\{\jn>j^*\}}\right)
&\leq C''' \tilde{B}^{2d/(2s+d)}n^{-2s/(2s+d)},
\end{align}
since $\tilde{B}>1.$\\
Finally by  \eqref{KoncowaDlaI}, \eqref{KoniecDlaII} and \eqref{eq:DwieCzesci} we obtain the Theorem.

\subsection{Proof of Lemma \ref{krotko}} \label{A8}
The idea of proof is from \cite{Hull}. Note that
\begin{align*}
&\mathbb{E} \ll{f_n(j) - \mathbb{E}f_n(j)}^4 =\\
&= \mathbb{E}\left( \ints \left|f_n(j)(x)-\mathbb{E}f_n(j)(x)\right|^2 \dS{x}  \right)^2\\
& = \mathbb{E}\left( \ints\ints\left| \frac 1n\sum_{i=1}^n (K_j(x,X_i) - K_jf(x))\right|^2  
\left| \frac 1n\sum_{l=1}^n (K_j(y,X_l) - K_jf(y))\right|^2  \dS{x}\dS{y}\right)^2.
\end{align*}
For $x\in\Sd$ we denote $Y_i(x) = K_j(x,X_i) - K_jf(x).$ Hence $\mathbb{E}Y_i(x) = 0.$
Thus
\begin{align*}
&\mathbb{E} \ll{f_n(j) - \mathbb{E}f_n(j)}^4 = \frac{1}{n^4} \mathbb{E}\left(\ints\ints \left|\sum_{i=1}^n Y_i(x)\right|^2 \left| \sum_{l=1}^n Y_l(y)\right|^2 \dS{x}\dS{y}\right)\\
&= \frac{1}{n^4} \mathbb{E}\Bigg(\ints\ints \Bigg(\sum_{i,l=1}^n Y_i^2(x)Y_l^2(y) 
    + 2 \sum_{i=1}^n\sum_{l<m} Y_i^2(x)Y_l(y)Y_m(y) +
    \\
&  \hspace{1.7cm} + 2 \sum_{l=1}^n\sum_{i<k} Y_l^2(y)Y_i(x)Y_k(x)+ 4 \sum_{i<k}\sum_{l<m} Y_i(x)Y_k(x)Y_l(y)Y_m(y)
   \Bigg)  \dS{x}\dS{y}\Bigg).
   \end{align*} 
   Put
\begin{equation}\label{eq:I_1234  proof:main}   
 \mathbb{E} \ll{f_n(j) - \mathbb{E}f_n(j)}^4   =\frac{1}{n^4}(I_1+I_2+I_3+I_4).
\end{equation}
Now we will estimate all components of  \eqref{eq:I_1234  proof:main}.
We will use the following inequality which is a consequence 
of Jensen's inequality for all $x\in\Sd$
\begin{equation}\label{eq:Jensen}
(K_jf(x))^2=\left( \ints K_j(x,u) f(u)\dS{u}\right)^2 \leq \ints K_j^2(x,u) f(u)\dS{u} = \mathbb{E} K_j^2(x,X).
\end{equation}
 We start with $I_1.$
\begin{align*}
I_1&=\mathbb{E}\left(\ints\ints \sum_{i,l=1}^n Y_i^2(x)Y_l^2(y) \dS{x}\dS{y} \right) \\
&=\sum_{i=1}^n \mathbb{E}\left(\ints\ints Y_i^2(x)Y_i^2(y) \dS{x}\dS{y}\right) 
+\sum_{1\leq i\neq l \leq n}^n \hspace{-0.3cm}\mathbb{E}\left(\ints\ints Y_i^2(x)Y_l^2(y)\dS{x}\dS{y}\right) 
\end{align*}
Put
\begin{equation}\label{eq:I_1_przed proof:main}
I_1=\sum_{i=1}^n I_1^{i,i} +  \sum_{1\leq i\neq l \leq n} I_1^{i,l}.
\end{equation}
For  $i=1, 2, \ldots, n$  using the inequality  $(a+b)^2\leq 2(a^2+b^2),$ and \eqref{eq:Jensen} we get for any $i,l$
\begin{align*}
I_1^{i,l}&=\mathbb{E}\left(\ints\ints Y_i^2(x)Y_l^2(y)\dS{x}\dS{y}\right)\nonumber\\
&= \mathbb{E}\left(\ints\ints            \Big[K_j(x,X_i) - K_jf(x)\Big] ^2   \Big[K_j(y,X_l) - K_jf(y)\Big] ^2           \dS{x}\dS{y}\right)\\
&\leq  4\mathbb{E}\left(\ints\ints            \Big[K_j^2(x,X_i) + (K_jf(x))^2\Big]   \Big[K_j^2(y,X_l) + (K_jf(y))^2\Big]           \dS{x}\dS{y}\right)\\
&\leq  4\ints\ints   \bigg\{ \mathbb{E}\Big[K_j^2(x,X_i)K_j^2(y,X_l)\Big]  +  3  \mathbb{E}\Big[K_j^2(x,X_i)\Big]\mathbb{E}\Big[K_j^2(y,X_l)\Big]    \bigg\}      \dS{x}\dS{y}.
\end{align*}
Now
\begin{align}\label{eq:I_1^i__2 proof:main}
I&_1^{i,i}\leq 4\ints\ints \ints K_j^2(x,u)K_j^2(y,u)f(u) \dS{x}\dS{y}\dS{u}+\nonumber\\
& \hspace{2cm} +  12 \ints  \mathbb{E}\Big[K_j^2(x,X_i)\Big] \dS{x} \ints \mathbb{E}\Big[K_j^2(y,X_i)\Big]\dS{y}.
\end{align}
For all $x\in\Sd$ from Lemma \ref{lem:zalozenie}
\begin{align}\label{eq:intsEKj^2 proof:main}
\ints  \mathbb{E}\Big[K_j^2(x,X_i)\Big] \dS{x} &= \ints \ints K_j^2(x,u) f(u) \dS{x}\dS{u}\nonumber\\ 
&\leq \ints f(u) \cdot \D2^{jd}\dS{u}= \D2^{jd}.
\end{align}
Hence using \eqref{eq:I_1^i__2 proof:main} and Lemma \ref{lem:zalozenie} we get
\begin{align}\label{eq:I_1^i__caly proof:main}
I_1^{i,i}&\leq  4\ints\ints \ints K_j^2(x,u)K_j^2(y,u)f(u) \dS{x}\dS{y}\dS{u} + 12 (\D2^{jd})^2\nonumber\\
&\leq 4\ints f(u) \left(\D2^{jd}\right)^2\dS{u} + 12\D^22^{2jd}=16 \D^2 2^{2jd}.
\end{align}
For all $1\leq i\neq l \leq n$ random variable $K_j^2(x,X_i)$ and $K_j^2(y,X_l)$ are independent using  \eqref{eq:intsEKj^2 proof:main}, we obtain
\begin{align}\label{eq:I_1^i,l__caly proof:main}
I_1^{i,l}&=16\ints   \mathbb{E}\Big[K_j^2(x,X_i)\Big]\dS{x} \ints \mathbb{E}\Big[K_j^2(y,X_l)\Big] \dS{y} \leq 16 (\D2^{jd})^2.
\end{align}
Consequently  \eqref{eq:I_1_przed proof:main} and \eqref{eq:I_1^i__caly proof:main} with  \eqref{eq:I_1^i,l__caly proof:main} give
\begin{align}\label{eq:I_1_caly proof:main}
I_1=\sum_{i=1}^n I_1^{i,i} +  \sum_{1\leq i\neq l \leq n}^n I_1^{i,l}\leq\sum_{i=1}^n 16 \D^2 2^{2jd} +  \sum_{1\leq i\neq l \leq n} 16 \D^2 2^{2jd}= 16 n^2  \D^2 2^{2jd}.
\end{align}
Let us estimate $I_2$ and $I_3.$
Since $Y_1, \ldots, Y_n$ are independent, we get
\begin{align}\label{eq:I_2__caly proof:main}
I_2 = 2 \mathbb{E}\left(\ints\ints\sum_{i=1}^n\sum_{l<m} Y_i^2(x)Y_l(y)Y_m(y)  \dS{x}\dS{y}\right)=0.
\end{align}
Similar
\begin{align}\label{eq:I_3__caly proof:main}
I_3 = 2 \mathbb{E}\left(\ints\ints\sum_{l=1}^n\sum_{i<k} Y_l^2(y)Y_i(x)Y_k(x)  \dS{x}\dS{y}\right) =0.\\
&\nonumber
\end{align}
To finish the proof of Lemma we need to estimate  $I_4.$,
\begin{align}\label{eq:I_4__1 proof:main}
I_4= 4 \sum_{i<k}\sum_{l<m} \mathbb{E}\left(\ints\ints   Y_i(x)Y_k(x)Y_l(y)Y_m(y)\dS{x}\dS{y}\right)= 4 \sum_{i<k}\sum_{l<m} I_4^{i,k,l,m}.
\end{align}
Note that for $i<k$, $l<m$ and assuming $i=l$ i $k=m$ we get
\begin{align}\label{eq:I_4__2 proof:main}
I_4^{i,k,l,m}= \ints\ints  \Big\{\mathbb{E}\big[ Y_i(x)Y_i(y)\big]\Big\}^2\dS{x}\dS{y}.
\end{align}
From Jensen's inequality
\eqref{eq:I_1^i__caly proof:main}
 we get
\begin{align}\label{eq:I_4__5 proof:main}
I&_4^{i,k,l,m}\leq I_1^{i,i}\leq 16 \D^22^{2jd}.
\end{align}
If  $i<k$ and $l<m$ assuming that $i\neq l$ or $k\neq m$ we get 
\begin{align}\label{eq:I_4__6 proof:main}
I_4^{i,k,l,m}=\ints\ints   \mathbb{E}\big[Y_i(x)Y_k(x)Y_l(y)Y_m(y)\big]\dS{x}\dS{y}=0.
\end{align}
Now by \eqref{eq:I_4__5 proof:main}, \eqref{eq:I_4__6 proof:main} and \eqref{eq:I_4__1 proof:main} we get
\begin{align}\label{eq:I_4__caly proof:main}
I_4&= 4 \sum_{i<k}\sum_{l<m} I_4^{i,k,l,m}\leq n^2\cdot 16 \D^22^{2jd} = 16 n^2\D^22^{2jd}.
\end{align}
By \eqref{eq:I_1_caly proof:main}, \eqref{eq:I_2__caly proof:main}, \eqref{eq:I_3__caly proof:main} and \eqref{eq:I_4__caly proof:main} applying  \eqref{eq:I_1234  proof:main} we get
\begin{align*}
\mathbb{E} \ll{f_n(j) - \mathbb{E}f_n(j)}^4 &= \frac{1}{n^4}(I_1+I_2+I_3+I_4)\leq \frac{32}{n^4}n^2 \D^22^{2jd}\nonumber\\
&= 32\D^2\left(\frac{2^{jd}}{n}\right)^2 = 32\D^2 \sigma^4(j,n).
\end{align*}
Finally to finish the proof we use above inequality and \eqref{def:bias} 
with \eqref{eq:E()^4 proof:main}. Indeed
\[
\mathbb{E} \ll{f_n(j) - f}^4 \leq  4\left(32\D^2 \sigma^4(j,n) + (\sqrt{2}\cdot 2^{-js} \nb{f})^4\right).
\]

\vskip 14pt
\noindent {\it Acknowledgements.}
The work of B. \'Cmiel was partially supported by the Faculty of Applied Mathematics AGH UST dean grant for PhD students and young researchers within subsidy of Ministry of Science and Higher Education.
\par


\bibliographystyle{amsplain}

\end{document}